\documentclass[a4paper,12pt]{article}
\usepackage{amssymb,amsmath,graphics,pstricks}
\newtheorem{thm}{Theorem}[section]
\newenvironment{dwd}{\par\noindent{\bf Proof.}}{\par\rightline{$\blacksquare$}}

\newtheorem{theo}[thm]{Theorem}
\newtheorem{pros}[thm]{Proposition}  
\newtheorem{coro}[thm]{Corollary}

\newtheorem{lema}[thm]{Lemma}
\newtheorem{defi}[thm]{Definition}
\newtheorem{rema}[thm]{Remark}

\def\be#1\ee{\begin{equation}#1\end{equation}}
\newcommand{\ba}{\begin{eqnarray} }
\newcommand{\ea}{\end{eqnarray} }
\def\bt#1\et{\begin{theo}#1\end{theo}}
\def\bl#1\el{\begin{lema}#1\end{lema}}
\def\bp#1\ep{\begin{pros}#1\end{pros}}
\def\bd#1\ed{\begin{defi}#1\end{defi}}

\def\ccB{{\cal B}}

\def\ccF{{\cal F}}

\def\ccN{{\cal N}}

\def\ccS{{\cal S}}

\def\va{\varepsilon}
\def\ra{\rightarrow}

\def\C{{\mathbb C}}
\def\E{\mathbf{E}}
\def\P{\mathbf{P}}

\def\N{{\mathbb N}}

\def\R{{\mathbb R}}
\def\Z{{\mathbb Z}}

\def\ls{\leqslant}
\def\gs{\geqslant}
\def\for{\mbox{for}}

\numberwithin{equation}{section}

\begin{document}

\title{\bf The Kendall's Theorem and its Application to the Geometric Ergodicity of Markov Chains }

\author{Witold Bednorz}
\date{}

\maketitle

\renewcommand{\thefootnote}{}
\footnote{Department of Mathematic, University of Warsaw, Banacha 2, 02-097 Warsaw, Poland}
\footnote{Research partially supported by  MNiSW Grant 
N N201 387234 }
\footnote{The paper was prepared during stay in IM PAN}
\footnote{2010 \emph{Mathematics Subject Classification}: Primary 60J20; Secondary 60K05;65C05.}
\footnote{\emph{Geometric Ergodicity; Renewal Theory; Markov Chain Monte Carlo}}
\renewcommand{\thefootnote}{\arabic{footnote}}
\setcounter{footnote}{0}

\begin{abstract}
In this paper we prove a sharp quantitative version of the Kendall's Theorem. The Kendal Theorem states that under 
some mild conditions imposed on a probability distribution on positive integers (i.e. probabilistic sequence) one can prove
convergence of its renewal sequence. Due to the well-known property - the first entrance last exit decomposition - such results are of interest 
in the stability theory of time homogeneous Markov chains. In particular the approach may be used to measure 
rates of convergence of geometrically ergodic Markov chains and consequently implies
estimates on convergence of MCMC estimators. 
\end{abstract}

\section{Introduction} \label{sect0}

Let $(X_n)_{n\gs 0}$ be a time-homogeneous Markov chain on a measurable space $(\ccS,\ccB)$, 
with transition probabilities $\P^n(x,\cdot)$, $n\gs 0$ and a unique stationary measure $\pi$. Let $P$ be the transition operator
given on the Banach space of bounded measurable functions on $(\ccS,\ccB)$ by
$P f(x)=\int f(y)\P(x,dy)$. Under mild conditions imposed on $(X_n)_{n\gs 0}$ 
the chain is ergodic, i.e.
\be\label{ergo}
\|\P^n(x,\cdot)-\pi(\cdot)\|_{TV}\ra 0,\;\;\mbox{as}\;\;n\ra\infty,
\ee 
for all starting points $x\in \ccS$ in the usual total variation norm 
$$
\|\mu\|_{TV}=\sup_{|f|\ls 1}|\int f d\mu|,
$$
where $\mu$ is a real measure on $(\ccS,\ccB)$.
It is known that the aperiodicity, the Harris recurrence property  and the finiteness of $\pi$ are equivalent to (\ref{ergo}), (see Theorem 13.0.1 in \cite{MT1}). Consequently 
the recurrence property is necessary to prove the convergence of $X_n$ distributions to the invariant measure in the total variation norm regardless of the starting point $X_0=x$.
In applications (see \cite{Num1}) there is required a stronger form of the result, namely we expect the exponential rate of the convergence and a reasonable method to estimate 
this quantity.
\smallskip

\noindent
One of the possible generalizations of the total variation convergence is considering functions controlled from above by $V:\ccS\ra \R$, $V\gs 1$, $\pi(V)<\infty$ therefore we 
refer to $B_V$ as the Banach space of measurable functions on $(\ccS,\ccB)$, such that $\sup_{x\in \ccS}|f(x)|/V(x)<\infty$ with the norm 
$$
\|f\|_V:=\sup_{x\in \ccS}\frac{|f(x)|}{V(x)}.
$$
Then instead of the total variation distance one can use
$$
\|\mu\|_{V}:=\sup_{|f|\ls V}|\int f d\mu|.
$$
The geometric convergence of $\P^n(x,\cdot)$ to a unique stationary measure $\pi$,
means there exists $\rho_V  <r\ls 1$ such that
\be\label{conv1}
\|(P^ng)(x)-\int g d\pi\|_V\ls M_V(r) r^n\|g\|_V\;\;g\in B_V,
\ee
where $\rho_V$ is the spectral radius of $(P-1\otimes \pi)$ acting on $(B_V,\|\cdot\|_B)$
and $M_V(r)$ is the optimal constant. In applications we often work with test functions $g$ from
a smaller space $B_W$, where $W:\ccS\ra \R$ and $1\ls W\ls V$. In this case we expect
$$
\|(P^ng)(x)-\int g d\pi\|_V\ls M_W(r) r^n \|g\|_W,\;\;g\in B_W,
$$
which is valid at least on $\rho_V\ls r\ls 1$, and $M_W(r)$ is the optimal constant. 
The most important case is when $W\equiv 1$, i.e. we consider not necessarily uniform geometric
convergence in the total variation norm.
\smallskip

\noindent
Whenever it exists we call $\rho_V$ the convergence rate of geometric ergodicity for the chain $(X_n)_{n\gs 0}$. For a class of examples one can prove the geometric convergence (see Chapter 15 in \cite{MT1}) and 
it is closely related to the existence of the exponential moment of the return time for a set $C\in\ccB$ of positive $\pi$-measure. 
\smallskip

\noindent The main tool to measure the convergence rate of the geometric ergodicity
is the drift condition, i.e. the existence of Lyapunov function $V:\ccS\ra \R$, $V\gs 1$, which is contracted outside a small set $C$.
The standard formulation of the required properties is the following:
\begin{enumerate}
\item \label{A1}\textit{Minorization condition.} There exist $C\in \ccB$, $b,\bar{b}>0$ and a probability measure $\nu$ on $(\ccS,\ccB)$ such that
$$
\P(x,A)\gs \bar{b}\nu(A)
$$ 
for all $x\in C$ and $A\in \ccB$.
\item\label{A2} \textit{Drift condition.} 
There exist a measurable function $V:S\ra[1,\infty)$ and constants
$\lambda<1$ and $K<\infty$ satisfying
$$
PV(x)\ls 
\left\{\begin{array}{ll}
\lambda V(x) & \mbox{if}\; x\not\in C\\
K            & \mbox{if}\; x\in C.    
\end{array}\right.
$$
\item\label{A3} \textit{Strong aperiodicity} $\bar{b}\nu(C)\gs b>0$ 
\end{enumerate}
The first property means there exists a small set $C$ on which the regeneration of
$(X_n)_{n\gs 0}$ takes place (see Chapter 5 in \cite{MT1}). The assumption is relatively week since each Harris recurrent chain admits the existence of a small set at least for some of its $m$-skeletons  (i.e. processes $(X_{nm})_{n\gs 0}$, $m\gs 1$) - see Theorem 5.3.2 in \cite{MT1}. The small set existence is used in the split chain construction (see Section \ref{sect4} and cf. \cite{Num} for details) to 
extend $(X_n)_{n\gs 0}$ to a new Markov Chain on a larger probability space $\ccS\times \{0,1\}$,
so that $(C,1)$ is a true atom of the new chain and its marginal distribution on $(\ccS,0)$ equals the distribution of $(X_n)_{n\gs 0}$.
The second condition reads as the existence of a Lyapunov function $V$ which is contracted by the semigroup related operator $P$ 
with the rate $\lambda<1$, for all points outside the small set.
Finally the strong aperiodicity means that the regeneration set $C$ is of positive measure  
for the basic transition probability, so the regeneration can take place in one turn.
\smallskip

\noindent
Our main result concerns ergodic Markov chains. Since the approach is based on the reduction to the study of a renewal sequence, we
first consider the atomic case and show how the idea works in this special case. Then we turn to the general case where the split chain construction is 
required. The results of the type are used whenever exact estimates on the ergodicity are required cf. \cite{Ad-Bed}, \cite{Niem1} and \cite{Niem2}.

\smallskip 

\noindent
The organization of the paper is as follows: in Section \ref{sect1a} we discuss the atomic case where we describe our main results in this case, then
in Section \ref{sect1b} we say what can be done with a method of chain splitting in the general case.
In  Section \ref{sect2} we 
give the proof of our main result - Theorem \ref{kendal1}; in Section \ref{sect3} we discuss how the result improves the previously known estimation methods. 
In section \ref{sect4} we show how the Kendal type results can be used in Markov theory to estimate rates of convergence. We give a short argument for both the atomic and non atomic case, leaving the tedious computation of some estimates of constants (which improves what has been known) to the Appendix A. Then in Appendix Btes of constants (which improves what has been known) to the Appendix A. Then in Appendix B
we analyze the result for typical toy examples.

\section{The atomic case}\label{sect1a}

\noindent
For this section we assume that $\bar{b}=1$. Note that in this setting one can rewrite the minoriziation condition \ref{A1} as  
$$
\P(x,A)=\nu(A),\;\;\mbox{for all}\; x\in C.
$$ 
which implies that $C$ is an atom and $\nu=\P(a,\cdot)$, for any $a\in C$. It remains to 
translate conditions \ref{A2}-\ref{A3} into a simpler form which can be used later to prove the geometric ergodicity.
Let $\tau=\tau(C)=\inf\{n\gs 1:\;X_n\in C\}$ and
$u_n=\P_a(X_n\in C)$, for $n\gs 0$. Then $u_n$ is the renewal sequence that corresponds to the increment sequence $b_n=\P_a(\tau=n)$
for $n\gs 1$. Note that in particular whenever we expect ergodicity $\lim_{n\ra\infty}u_n$ exists and is equal $u_{\infty}=\pi(C)$. Following \cite{Bax} we define function $G(r,x)=\E_x r^{\tau}$, for all $x\in S$ and $0<r\ls \lambda^{-1}$.
The main property of $G(r,x)$ is that it is the lower bound for $V(x)$ on the set $\ccS\backslash C$, namely
we have that (cf. Proposition 4.1 in \cite{Bax})
\begin{pros}\label{pro1} Assume only drift condition (\ref{A2}).
\begin{enumerate}
\item For all $x\in S$, $\P_x(\tau<\infty)=1$.
\item For $1\ls r\ls \lambda^{-1}$
$$
G(r,x)\ls\left\{
\begin{array}{ll}
V(x) & \mbox{if}\;x\not \in C,\\
rK & \mbox{if}\;x\in C.
\end{array}\right.
$$
\end{enumerate}
\end{pros}
The renewal approach is based on the first entrance last exit property. To state it we 
need additional notation $H_W(r,x)=\E_x(\sum^{\tau}_{n=1}r^nW(X_n))$, for $r>0$ for which the definition makes sense.
We have that (cf. Proposition 4.2 in \cite{Bax}) 
\begin{pros}\label{pro2}
Assume only that the Markov chain is geometrically ergodic with (unique) invariant probability measure $\pi$, that $C$ is an atom and
that $W:\ccS\ra\R$ is such that $W\gs 1$. Suppose $g:S\ra\R$ satisfies $\|g\|_W\ls 1$, then for all $r\gs 1$ 
for which right-hand sides are finite:
\[
\begin{split}
&\sup_{|z|= r}|\sum^{\infty}_{n=1}(P^ng(a)-\int g d\pi)z^n|\ls\\
&\ls H_W(r,a)\sup_{|z|\ls r}|\sum^{\infty}_{n=0}(u_n-u_{\infty})z^n|+\pi(C)\frac{H_W(r,a)-rH_W(1,a)}{r-1},
\end{split}
\]
for all $a\in C$ and
\[
\begin{split}
&\sup_{|z|= r}|\sum^{\infty}_{n=1}(P^ng(x)-\int g d\pi)z^n|\ls\\
&\ls H_W(r,x)+G(r,x)H_W(r,a)|\sup_{|z|\ls r}\sum^{\infty}_{n=0}(u_n-u_{\infty})z^n|+\\
&+\pi(C)\frac{H_W(r,a)-rH_W(1,a)}{r-1}G(r,x)+\pi(C)H_W(1,a)\frac{r(G(r,x)-1)}{r-1}, 
\end{split}
\]
for all $x\not\in C$. 
\end{pros}
Now the problem splits into two parts: in the first one we have to provide some estimate on $H_{W}(r,x)$, $x\in \ccS$ on the interval $1\ls r\ls \lambda^{-1}$,
and it is of meaning when we want to obtain a reasonable bound on $M_W(r)$, whereas in the second we search for $r_0$ a lower bound for the inverse of the 
radius of convergence of $\sum^{\infty}_{n=0}(u_n-u_{\infty})z^n$, and then for some upper bound $K_0(r)$ on $\sup_{|z|= r}|\sum^{\infty}_{n=0}(u_n-u_{\infty})z^n|$, for $r< r_0$.
\smallskip
 
\noindent 
As for the first issue we acknowledge two cases.  The simplest situation is when $W\equiv 1$ and therefore $H_{1}(r,x)=r(G(r,x)-1)/(r-1)$, $H_1(1,a)=\E_a \tau=\pi(C)^{-1}$, 
which allows to slightly improve estimates on $H_1(r,x)$ (cf. Proposition 4.2 in \cite{Bax}).
\begin{pros}\label{pro3}
Assume only drift condition (\ref{A2}).
\begin{enumerate}
\item For $1\ls r \ls \lambda^{-1}$
$$
H_1(r,x) \ls\left\{
\begin{array}{ll}
\frac{r\lambda(V(x)-1)}{1-\lambda} & \mbox{if}\;x\not \in C,\\
\frac{r(K-\lambda)}{1-\lambda} & \mbox{if}\;x\in C.
\end{array}\right.
$$
\item and for $1\ls r\ls \lambda^{-1}$
$$
\frac{H_1(r,a)-rH_1(1,a)}{r-1}\ls \frac{r\lambda(K-1)}{(1-\lambda)^2}.
$$
\end{enumerate}
\end{pros}
Combining estimates from Propositions \ref{pro1}
and \ref{pro3} with Proposition \ref{pro2} we obtain
\bt\label{thm1}
Suppose $(X_n)_{n\gs 0}$ satisfies conditions \ref{A1}-\ref{A3} with $\bar{b}=1$. Then $(X_n)_{n\gs 0}$ is geometrically ergodic - it verifies (\ref{conv1}) and we have 
the following bounds on $\rho_V$, $M_1$:
\begin{align*}
& \rho_V\ls r_0^{-1}\\
& M_1(r)\ls \frac{2r\lambda}{1-\lambda}+\frac{r\lambda(K-1)}{(1-\lambda)^2}+\frac{r(K-\lambda)}{1-\lambda} K_0(r),
\end{align*}
where $r_0=r_0(b,\lambda^{-1},\lambda^{-1}K)$ and $K_0(r)=K_0(r,b,\lambda^{-1},\lambda^{-1}K)$ are defined in Corollaries \ref{cori1},\ref{cori2}.
\et
On the other hand when $W\equiv V$ we have weaker bounds on $H_V(r)$, which are given in Proposition 4.2
in \cite{Bax}.
\begin{pros}\label{pro4} Assume only drift condition (\ref{A2}).
\begin{enumerate}
\item For $1\ls r \ls \lambda^{-1}$
$$
H_V(r,x) \ls\left\{
\begin{array}{ll}
\frac{r\lambda (V(x)-1)}{1-r\lambda} & \mbox{if}\;x\not \in C,\\
\frac{r(K-r\lambda)}{1-r\lambda} & \mbox{if}\;x\in C.
\end{array}\right.
$$
in particular $H_V(1,x)\ls \frac{K-\lambda}{1-\lambda}$ for all $x\in C$.
\item and for $1\ls r\ls \lambda^{-1}$
$$
\frac{H_V(r,a)-rH_V(1,a)}{r-1}\ls \frac{r\lambda  (K-1) }{(1-\lambda)(1-r\lambda)}.
$$
\end{enumerate}
\end{pros}
Applying Proposition \ref{pro4} instead of \ref{pro3} we obtain a similar result to Theorem \ref{thm1},
yet with a worse control on $M_W(r)$ (that necessarily goes to infinity near $r=\lambda^{-1}$).
\bt\label{thm2}
Suppose that $(X_n)_{n\gs 0}$ satisfies conditions \ref{A1}-\ref{A3} with $\bar{b}=1$. Then $(X_n)_{n\gs 0}$ is geometrically ergodic - it verifies (\ref{conv1}) and we have 
the following bounds on $\rho_V$, $M_V$:
\begin{align*}
& \rho_V\ls r_0^{-1}\\
& M_V(r)\ls \frac{r\lambda}{1-r\lambda}+
\frac{r\lambda(K-\lambda)}{(1-\lambda)^2}+\frac{r\lambda(K-1)}{(1-\lambda)(1-r\lambda)}+\frac{r(K-r\lambda)}{1-r\lambda}K_0(r),
\end{align*}
where $r_0=r_0(b,\lambda^{-1},\lambda^{-1}K)$ and $K_0(r)=K_0(r,b,\lambda^{-1},\lambda^{-1}K)$ are defined Corollaries \ref{cori1},\ref{cori2}.
\et
The second part concerns the study of a renewal process. As we have noted whenever the condition (\ref{A2}) holds we can handle with all quantities in Proposition \ref{pro2} but $\sup_{|z|\ls r}|\sum^{\infty}_{n=0}(u_n-u_{\infty})z^n|$.
Let $(\tau_k)_{k\gs 0}$ denote subsequent visits of the Markov chain to $C$, where we assume that $\tau_0=0$ (so the chain starts from $C$).
The renewal process is defined by $V_m=\inf\{\tau_n-m :\; \tau_n\gs m\}$, $m\gs 0$.
Clearly $\tau_k-\tau_{k-1}$, $k\gs 1$ forms a family of independent random variables of the same distribution 
$\P_a(\tau_k-\tau_{k-1}=n)=\P_a(\tau=n)=b_n$, $n\gs 1$. By the definition we have the equality $\P(V_n=0)=u_n$, i.e. the probability that the process renew in time $n$
equals $u_n$. Observe that 
$u_0=1$ and $u_{n}=\sum^n_{k=1}u_{n-k}b_{k}$, hence denoting 
$b(z)=\sum^{\infty}_{n=1}b_nz^n$, $u(z)=\sum^{\infty}_{n=0}u_nz^n$, for $z\in\C$, we can state the renewal equation in the
following form
\be\label{renewal}
u(z)=1/(1-b(z)),\;\;\for\; |z|<1.
\ee
The equation means we can study the convergence of $u_n$ to $u_{\infty}$ in terms of properties of $(b_n)_{n\gs 1}$.
Note that $b_1=\P_a(\tau=1)\gs b$ and $b(\lambda^{-1})\ls \lambda^{-1} K$ by the conditions respectively \ref{A1} and \ref{A2}. 
Historically, the first result that matches these properties with the geometric ergodicity was due to Kendall \cite{Ken} who proved that:
\bt\label{kendal}
Assume that $b_1>0$ and $\sum^{\infty}_{n=1}b_n r^n<\infty$ for some $r>1$. Then the limit 
$u_{\infty}=\lim_{n\ra\infty}u_n$ exists and is equal $u_{\infty}=(\sum^{\infty}_{n=1}nb_n)^{-1}$,
moreover the radius of convergence of $\sum^{\infty}_{n=0}(u_n-u_{\infty})z^n$ is strictly greater than $1$. 
\et
Although the result shows that drift condition implies the geometric ergodicity of an aperiodic Markov Chain, 
it is not satisfactory in the sense that it does not provide neither estimates on the radius of convergence nor a deviation inequality which one could use
on the disc where the convergence holds. The Kendall's theorem was improved first in \cite{MT2} and then in \cite{Bax} (see Theorem 3.2). 
There are also several results where some additional assumptions on the distribution of $\tau$ are studied. For example in \cite{BeLu} there is described how to provide an optimal bound on the rate of convergence, which gives some computable bound on the value, yet under additional conditions on the $\tau$ distribution. Whenever the general Kendall's problem is considered the bounds obtained in the mentioned papers  
are still far from being optimal or easy to use. The goal of the paper is to give a sharp estimate on
the rate of convergence which significantly improves on the previous results. Our approach is based on 
introducing $\pi(C)$ as a parameter, namely we prove that the following result holds:  
\bt\label{kendal1}
Suppose that $(b_n)^{\infty}_{n\gs 1}$ verifies $b_1\gs b >0$, $b(r)=\sum^{\infty}_{n=1}b_n r^n<\infty$, for some $r>1$.
Then $u_{\infty}=(\sum^{\infty}_{n=1}n b_n)^{-1}$ and 
$$
\sup_{|z|=r}|\sum^{\infty}_{n=0}(u_n-u_{\infty})z^n|\ls \frac{c(r)-c(1)}{c(1)(r-1)([(1-b)D(\alpha)-c(r)+c(1)]_{+})},
$$
where $c(r)=\frac{b(r)-1}{r-1}$, $c(1)=u_{\infty}^{-1}=\pi(C)^{-1}$ and 
$$
D(\alpha)=\frac{|1+\frac{b}{1-b}(1-e^{\frac{i\pi}{1+\alpha}})|-1}{|1-e^{\frac{i\pi}{1+\alpha}}|},\;\;\mbox{where}\;\;\alpha=\frac{c(1)-1}{1-b}, 
$$
\et
Consequently whenever one can control $c(r)=(b(r)-1)/(r-1)$ from above, there is a bound on the rate of convergence for 
the renewal process. The simplest exposition is when $c(1)=\pi(C)^{-1}$ is known and we can control $c(r)$
in a certain point, i.e. $c(R)\ls N<\infty$, for some $R>1$. Observe that if $b(R)\ls L$, then due to $c(R)=\frac{b(R)-1}{R-1}$
we deduce that $c(R)\ls N=\frac{L-1}{R-1}$, which is our basic setting.
Note that by the H\"{o}lder inequality, for all $1\ls r\ls R$
$$
c(r)-c(1)\ls (c(1)-1)(\frac{c(r)-1}{c(1)-1}-1)\ls (1-b)\alpha(r^{\kappa(\alpha)}-1),
$$
where $\kappa(\alpha)=\log(\frac{N-1}{c(1)-1})/\log R=\log(\frac{N-1}{(1-b)\alpha})/\log R$,
$\alpha=(c(1)-1)/(1-b)$.  
\begin{coro}\label{cori1}
Suppose that $c(1)=\pi(C)^{-1}$ is given, $b_1\gs b$ and $b(R)\ls L$, then
\be\label{ww1}
r_0=\min\{R,(1+\frac{D(\alpha)}{\alpha})^{\frac{1}{\kappa(\alpha)}}\}. 
\ee
Moreover for $r<r_0$
$$
\sup_{|z|=r}|\sum^{\infty}_{n=0}(u_n-u_{\infty})z^n|\ls K_0(r)= \frac{\pi(C)(r^{\kappa(\alpha)}-1)}{(r-1)(\alpha^{-1}D(\alpha)-r^{\kappa(\alpha)}+1)},
$$
\end{coro}
\begin{rema}\label{rasta1}
Observe that the bound $(1+\frac{D(\alpha)}{\alpha})^{\frac{1}{\kappa(\alpha)}}$ is increasing 
with $b$ assuming that $L,R,c(1)$ are fixed.
\end{rema}
In applications we have to treat $c(1)=\pi(C)^{-1}$ as a parameter.
The advantage of the approach is that there is a sharp upper bound on $c(1)$ or rather $\alpha=(c(1)-1)/(1-b)$. Using the inequality
\be\label{ola}      
R^{\alpha}=R^{(\sum^{\infty}_{n=1}(n-1)b_n)/(1-b)}\ls \frac{\sum^{\infty}_{n=2}b_nR^{n-1}}{1-b}\ls \frac{b(R)-bR}{(1-b)R}\ls \frac{L-bR}{(1-b)R},
\ee
we obtain that $\alpha \ls \alpha_0$, where $\alpha_0=\log(\frac{L-bR}{(1-b)R})/\log R$.
On the other hand if $b=b_1$, then $c(1)-1\gs 1-b$ and therefore due to 
Remark \ref{rasta1} we can always require that $c(1)-1\gs 1-b$.
Consequently to find an estimate on the rate of convergence we search 
$(1+\frac{D(\alpha)}{\alpha})^{\frac{1}{\kappa(\alpha)}}$, $\alpha\in [1,\alpha_0]$ for the
possible minimum. 
\begin{coro}\label{cori2}
Suppose that $b_1\gs b$ and $b(R)\ls L$. Then
\be\label{ww2}
r_0=\min\{R,\min_{1\ls \alpha\ls \alpha_0}(1+\frac{D(\alpha)}{\alpha})^{\frac{1}{\kappa(\alpha)}}\}.
\ee
Moreover for $r<r_0$
\be\label{ww3}
\sup_{|z|=r}|\sum^{\infty}_{n=0}(u_n-u_{\infty})z^n|\ls K_0(r)=\max_{1\ls \alpha \ls \alpha_0}\frac{r^{\kappa(\alpha)}-1}{(r-1)(\alpha^{-1}D(\alpha)-r^{\kappa(\alpha)}+1)}
\ee
\end{coro}
The above Corollary should be compared with Theorem 3.2. in \cite{Bax}. In section \ref{sect3} we prove that our result is always better than the previous one, moreover
we discuss the true benefit of the approach studying the limiting case, where $R\ra 1$.

\section{Non atomic case}\label{sect1b}

In the general case we assume $\bar{b}\ls 1$, which means that true atom may not exists. Therefore we have to use the split chain construction to create an atom on the extended probability space. In the Section \ref{sect4} we prove Theorems \ref{thm3} and \ref{thm4} that are equivalents of Theorems \ref{thm1} and \ref{thm2} in the case of general ergodic Markov chain.
Consequently applying the first entrance last exit decomposition for the split chain we reduce the question of the rate of convergence to the study of the 
renewal sequence for $(\bar{b}_n)_{n\gs 1}$ - the probability distribution of the return time to the artificial atom.
\smallskip

\noindent
Let $(\bar{u}_n)_{n\gs 0}$ be the corresponding renewal sequence for $(\bar{b}_n)_{n\gs 1}$. In the same way as in the atomic case let $\bar{b}(z)$, $\bar{u}(z)$, $z\in C$ be corresponding generating functions and 
$\bar{c}(z)=(\bar{b}(z)-1)/(z-1)$. Clearly
$\bar{b}_1=\bar{b}\nu(C)\gs b$, and $\bar{c}(1)=\bar{b}^{-1}\pi(C)$ so as in the atomic case we have a control on the limiting behavior of 
$\bar{c}(z)-\bar{c}(1)$, namely
applying Theorem \ref{kendal1} we obtain that whenever $\bar{c}(r)<\infty$, then
\be\label{jasmine0}
\sup_{|z|=r}|\sum^{\infty}_{n=0}(\bar{u}_n-\bar{u}_{\infty})z^n|\ls \frac{\bar{c}(r)-\bar{c}(1)}{\bar{c}(1)(r-1)([(1-b)D(\bar{\alpha})-\bar{c}(r)+\bar{c}(1)]_{+})},
\ee
where $\bar{c}(r)=\frac{\bar{b}(r)-1}{r-1}$, $\bar{c}(1)=\bar{u}_{\infty}^{-1}=\bar{b}^{-1}\pi(C)^{-1}$ and 
$$
D(\bar{\alpha})=\frac{|1+\frac{b}{1-b}(1-e^{\frac{i\pi}{1+\bar{\alpha}}})|-1}{|1-e^{\frac{i\pi}{1+\bar{\alpha}}}|},\;\;\mbox{where}\;\;\bar{\alpha}=\frac{\bar{c}(1)-1}{1-b}, 
$$
In this way the problem reduces to the estimate on $\bar{b}(r)$. The main difficulty is that in the non atomic case 
Theorem \ref{A2} provides only that for $R=\lambda^{-1}>1$
\be\label{jasmine1}
b_x(R)=\E_x R^{\tau}\ls L=KR,\;\;\mbox{for all}\; x\in C,
\ee
whereas we need a bound on the generic function of $(\bar{b}_n)_{n\gs 1}$.
We discuss the question in Section \ref{sect4}, showing in Proposition \ref{pro7} that 
for all $1\ls r\ls \min\{R,(1-b)^{-\frac{1}{1+\alpha_1}}\}$ the following inequality holds
\be\label{ing1}
\bar{b}(r)\ls L(r)=\max\{\frac{\bar{b}r}{1-(1-\bar{b})r^{1+\alpha_1}},\frac{br+(\bar{b}-b)r^{1+\alpha_2}}{1-(1-\bar{b})r}\},
\ee
where $\alpha_1=\log(\frac{L-\bar{b}R}{(1-\bar{b})R})/\log R$ and $\alpha_2=\log(\frac{L-(1-\bar{b}+b)R}{(\bar{b}-b)R})/\log R$. 
Moreover if $1+b\gs 2\bar{b}$ then simply
$$
L(r)=\frac{\bar{b}r}{1-(1-\bar{b})r^{1+\alpha_1}}.
$$
Using (\ref{ing1}) is the best what the renewal approach can offer to bound $\bar{b}(r)$. The meaning 
of the result is that there are only two generic functions that are important to bound $\bar{b}(r)$.
If $\bar{b}$ is close to $1$ then we are in the similar setting as in the atomic case and surely
one can expect the bound on $\bar{b}(r)$ of the form $\frac{br+(\bar{b}-b)r^{1+\alpha_2}}{1-(1-\bar{b})r}$, whereas
if $\bar{b}$ is far from $1$ only the split chain construction matters and the bound on $\bar{b}(r)$ should be like
$\frac{\bar{b}r}{1-(1-\bar{b})r^{1+\alpha_1}}$. 
\smallskip

\noindent
As in the atomic case we will need a bound on the $\bar{\alpha}=\frac{\bar{c}(1)-1}{1-b}$. We show in Corollary \ref{jasmine2}
that
\be\label{valbone1}
\bar{\alpha}\ls \bar{b}^{-1}\max\{\frac{1-\bar{b}}{1-b}(1+\alpha_1),\frac{1-\bar{b}}{1-b}+\frac{\bar{b}-b}{1-b}\alpha_2\}.
\ee
In fact the maximum equals $\bar{b}^{-1}\frac{1-\bar{b}}{1-b}(1+\alpha_1)$ if $1+b\gs 2\bar{b}$ and $\bar{b}^{-1}\frac{1-\bar{b}}{1-b}+\frac{\bar{b}-b}{1-b}\alpha_2$ otherwise.
\smallskip

\noindent
Now we turn to the basic idea for all the approach presented in the paper.
Observe that $\frac{\bar{c}(r)-1}{\bar{c}(1)-1}$ satisfies the H\"{o}lder inequality i.e. for $p+q=1$, $p,q> 0$
$$
(\frac{\bar{c}(r_1)-1}{\bar{c}(1)-1})^{p}(\frac{\bar{c}(r_2)-1}{\bar{c}(1)-1})^{q}\gs \frac{\bar{c}(r_1^pr_2^q)-1}{\bar{c}(1)-1},
$$
which means that $F_0(x)=\log(\frac{\bar{c}(e^x)-1}{\bar{c}(1)-1})$ is convex and $F_0(0)=0$. By (\ref{ing1}) we have that
$\bar{c}(e^x)\ls L(e^x)$ and hence
\be\label{antyma1}
F_0(x)\ls F_1(x)=\log(\frac{L(e^x)-e^x}{(1-b)\bar{\alpha}(e^x-1)}).
\ee
Therefore we can easily compute the largest possible function $\bar{F}(x)$ that satisfies the conditions:
\begin{enumerate}
\item $\bar{F}(x)\ls F_1(x)$ for $0\ls x \ls \min\{\log R,-\frac{1}{1+\alpha_1}\log(1-\bar{b})\}$; 
\item $\bar{F}(0)=0$ and $\bar{F}$ is convex; 
\item $\bar{F}$ is maximal over the functions with the properties 1-2, namely if there exists $F$ that satisfies
the above condition then $F(x)\ls \bar{F}(x)$ for all $0\ls x\ls \min\{\log R,-\frac{1}{1+\alpha_1}\log(1-\bar{b})\}$. 
\end{enumerate}
Let $x_0$ be the unique solution of the equation 
\be\label{antyma2}
F'_1(x)x=F_1(x).
\ee
Note that $x_0\ls -\frac{1}{1+\alpha_1}\log(1-\bar{b})$. If additionally $x_0\ls \log R$ then the optimal $\bar{F}(x)$
is of the form
\be\label{jasmine3}
\bar{F}(x)=\left\{\begin{array}{ll}
F'_1(x_0)x & \mbox{for all} \;\;0\ls x\ls x_0 \\
F_1(x) & \mbox{for all}\;\;x_0\ls x\ls \min\{\log R,-\frac{1}{1+\alpha_1}\log(1-\bar{b})\}
\end{array}\right.
\ee
otherwise if $x_0>\log R$ then 
\be\label{jasmine4}
\bar{F}(x)=\frac{F_1(\log R)}{\log R}x\;\; \mbox{for all}\;\; 0\ls x\ls \log R.
\ee
To make the notation similar to the atomic case let $\bar{\kappa}(\bar{\alpha},r)=\frac{\bar{F}(\log r)}{\log r}$. In particular if
$\log r <x_0\ls \log R$ then $\bar{\kappa}(\bar{\alpha},r)=F'_1(x_0)$ and similarly $\bar{\kappa}(\alpha,r)=\frac{\bar{F}(\log R)}{\log R}$ if
$\log R< x_0$. The above discussion leads to the following conclusion: 
\be\label{antyma4}
\bar{c}(r)-\bar{c}(1)\ls (1-b)\bar{\alpha} r^{\kappa(\bar{\alpha},r)},\;\;\mbox{for all}\;\;1\ls r\ls \min\{R,(1-\bar{b})^{-\frac{1}{1+\alpha_1}}\}
\ee
furthermore $\bar{\kappa}(\bar{\alpha},r)$ as a function of $r$ is constant at least on the part of the interval $[1,\min\{R,(1-\bar{b})^{-\frac{1}{1+\alpha_1}}\}]$.
Consequently applying (\ref{jasmine0} in the case where $\pi(C)$ is known we obtain our main estimate in the non-atomic case.
\begin{theo}\label{cori6}  
Suppose that $\bar{b}_1\gs b$ and $\bar{b}(r)$ satisfies (\ref{ing1}), and $\pi(C)$ is known.
Then
$$
\bar{r}_0=\min\{R,(1-\bar{b})^{-\frac{1}{1+\alpha_1}},\bar{r}_0(\alpha)\},
$$
where $\bar{r}_0(\alpha)$ is the unique solution of the equation
$$
r=(1+\frac{D(\bar{\alpha})}{\bar{\alpha}})^{\frac{1}{\bar{\kappa}(\bar{\alpha},r)}}.
$$
Moreover for $r<r_0$
$$
\sup_{|z|=r}|\sum^{\infty}_{n=0}(u_n-u_{\infty})z^n|\ls \bar{K}_0(r)= 
\frac{\bar{b}\pi(C)(r^{\bar{\kappa}(\bar{\alpha},r)}-1)}{(r-1)(\bar{\alpha} D(\bar{\alpha})-r^{\bar{\kappa}(\bar{\alpha},r)}+1)}.
$$
\end{theo}
\begin{rema}
Observe that if 
\be\label{antyma3}
\log(1+\frac{D(\bar{\alpha})}{\bar{\alpha}})/F_1'(x_0)\ls x_0\ls \log R,
\ee
then $\bar{r}_0=(1+\frac{D(\bar{\alpha})}{\bar{\alpha}})^{\frac{1}{F_1'(x_0)}}$.
Due to (\ref{antyma2}), the condition (\ref{antyma3}) is equivalent to $x_0\ls \log R$ and
$$
1+\frac{D(\bar{\alpha})}{\alpha}\ls \frac{L(e^{x_0})-e^{x_0}}{(1-b)\bar{\alpha}(e^{x_0}-1)}.
$$
\end{rema}
Therefore for a large class of examples we have a computable direct bound on the rate of convergence
even for general ergodic Markov chains. 
\smallskip

\noindent
If $\pi(C)$ is unknown then we have to treat it as a parameter and use a bound on $\bar{\alpha}$.
As for the upper bounds we can use (\ref{valbone1}), on the other hand
we show in Corollary \ref{jasmine2} that if $\bar{b}\nu(C)=b$ then $\bar{\alpha}\gs \bar{b}^{-1}$.
Since in the same way as in the atomic case 
$(1+\frac{D(\bar{\alpha})}{\bar{\alpha}})^{\frac{1}{\bar{\kappa}(\bar{\alpha},r)}}$
increases with $b$ assuming that $\bar{b},L,R$ are fixed, thus we can always assume $\bar{\alpha}\gs \bar{b}^{-1}$. Let $\bar{\alpha}_0=\max\{\frac{1-\bar{b}}{1-b}(1+\alpha_1),\frac{1-\bar{b}}{1-b}+\frac{\bar{b}-b}{1-b}\alpha_2\}$.
\begin{theo}\label{cori7}  
Let $\bar{b}_1\gs b$, and $\bar{b}(r)$ satisfies (\ref{ing1}). Then
$$
\bar{r}_0=\min\{R,(1-\bar{b})^{-\frac{1}{1+\alpha_1}}, \min_{\bar{b}^{-1}\ls \bar{\alpha}\ls \bar{b}^{-1}\bar{\alpha}_0}\bar{r}_0(\bar{\alpha})\},
$$
where $\bar{r}_0(\bar{\alpha})$ is the unique solution of the equation
$$
r=(1+\frac{D(\bar{\alpha})}{\bar{\alpha}})^{\frac{1}{\bar{\kappa}(\bar{\alpha},r)}}.
$$
Moreover for $r<r_0$          
$$
\sup_{|z|=r}|\sum^{\infty}_{n=0}(u_n-u_{\infty})z^n|\ls \bar{K}_0(r)= 
\max_{\bar{b}^{-1}\ls \bar{\alpha}\ls \bar{b}^{-1}\bar{\alpha}_0}
\frac{\bar{b}(r^{\bar{\kappa}(\bar{\alpha},r)}-1)}{(r-1)(\bar{\alpha}^{-1}D(\bar{\alpha})-r^{\bar{\kappa}(\bar{\alpha},r)}+1)}.
$$
\end{theo}
We show in examples that the approach presented in Theorems \ref{cori6}, \ref{cori7} is comparable with the coupling method (see Section 7 in  \cite{Bax} for short introduction).
Therefore we obtain the computable tool for the general question of rates of convergence of ergodic Markov chains under the geometric drift condition. 

\section{Proof of main result} \label{sect2}

In this section we give a proof of Theorem \ref{kendal1}.
\begin{dwd}[of Theorem \ref{kendal1}]
Let $b(z)$ and $u(z)$ be the complex generic functions for $b_i$, $i\gs 1$ and $u_i$, $i\gs 0$ sequences respectively. 
The main tool we use is the renewal equation (\ref{renewal}), i.e.
$$
1-b(z)=\frac{1}{u(z)},\;\;|z|<1.
$$
Note that the equation remains valid on the disc $|z|\ls R$ in the sense of analytic functions.
By Theorem \ref{kendal} we learn that $u_{\infty}<\infty$ and 
the renewal generic function $\sum^{\infty}_{n=0}(u_n-u_{\infty})z^n$ is convergent 
on some disc with radius greater than $1$. Denote
$c(z)=\frac{b(z)-1}{z-1}$ (cf. proof of Theorem 3.2 in \cite{Bax}) and observe that $c(z)$ is well defined on $|z|\ls R$, because $c(R)=\frac{b(R)-1}{R-1}=\frac{L-1}{R-1}<\infty$.
Since $u_{\infty}=c(1)^{-1}$ we have that
\begin{eqnarray}\label{main0}
&& \sum^{\infty}_{n=0}(u_n-u_{\infty})z^n=u(z)-\frac{1}{c(1)(1-z)}=\frac{1}{1-b(z)}-\frac{1}{c(1)(z-1)}=\nonumber\\
&&=\frac{1}{1-z}(\frac{1}{c(z)}-\frac{1}{c(1)})=\frac{c(z)-c(1)}{z-1}\frac{1}{c(1)c(z)}.
\end{eqnarray}
The main problem is to estimate $|c(z)|$ from below, to which goal we use the simple technique
\be\label{ineq1}
|c(re^{i\theta})|=|c(e^{i\theta})|-|c(re^{i\theta})-c(e^{i\theta})|=|c(e^{i\theta})|-c(r)+c(1).
\ee
Consequently the problem is reduced to the study of the lower bound on $|c(e^{i\theta})|$. We recall that by the definition $c_i=\sum_{j>i}b_j$ and 
$c(1)=\sum^{\infty}_{i=0}c_i$. To provide a sharp estimate in (\ref{ineq1}) we benefit from the fact that for $\frac{\pi}{l+1}<|\theta|\ls\frac{\pi}{l}$, $l\gs 1$,
there is a better control on the first $l$ summands in $c(e^{i\theta})=\sum^{\infty}_{j=1}c_je^{ij\theta}$.  First we note that  
$$
|c(e^{i\theta})|=\frac{|1-\sum^{\infty}_{j=1}b_je^{ij\theta}|}{|1-e^{i\theta}|}\gs
\frac{|1-\sum^l_{j=1}b_je^{ij\theta}|-\sum_{j>l}b_j}{|1-e^{i\theta}|},
$$
which is equivalent to
$$
|c(e^{i\theta})|\gs \frac{|c_l+\sum^l_{j=1}b_j(1-e^{ij\theta})|-c_l}{|1-e^{i\theta}|}.
$$
The geometrical observation gives that for $\frac{\pi}{l+1}<|\theta|\ls \frac{\pi}{l}$
$$
|c_l+\sum^l_{j=1}b_j(1-e^{ij\theta})|\gs |c_l+(\sum^l_{j=1}b_j)(1-e^{i\theta})|=
|c_l+(1-c_l)(1-e^{i\theta})|,
$$
hence we conclude that
$$
|c(e^{i\theta})|\gs c_l |1-e^{i\theta}|^{-1}(|1+(1-c_l)c_l^{-1}(1-e^{i\theta})|-1).
$$
Since $1-c_l\gs b$, for $l\gs 1$ we see that  
$$
|1+(1-c_l)c_l^{-1}(1-e^{i\theta})|\gs \sqrt{1+b c_l^{-2}|1-e^{i\theta}|^2}.
$$
It remains to verify that  $f(x)=x^{-1}[\sqrt{1+bx^2}-1]$ is  increasing,
which is assured by
\be\label{ineq2}
f'(x)=-x^{-2}(\sqrt{1+bx^2}-1)+x^{-2}\frac{bx^2}{\sqrt{1+4bx^2}}\gs 0.
\ee
Therefore we finally obtain that for $\frac{\pi}{l+1}<|\theta|\ls \frac{\pi}{l}$
\be\label{ccq}
|c(e^{i\theta})|\gs c_l|1-e^{\frac{i\pi}{l+1}}|(\sqrt{1+bc_l^{-2}|1-e^{\frac{i\pi}{l+1}}|^2}-1). 
\ee
Due to (\ref{ineq2}) and (\ref{ccq}), when estimating the global minimum of $|c(e^{i\theta})|$ it suffices  to find 
the bound from above on $c_l|1-e^{\frac{i\pi}{l+1}}|^{-1}$. We will show that
\be\label{ineq3}
c_l|1-e^{\frac{i\pi}{(l+1)}}|^{-1}\ls (1-b)|1-e^{\frac{i\pi}{1+\alpha}}|^{-1},
\ee
where we recall that $\alpha=(c(1)-1)/(1-b)$. First observe that (\ref{ineq3}) is trivial for $l\ls \alpha$,
since $c_l\ls (1-b)$ and $|1-e^{\frac{i\pi}{l+1}}|\gs |1-e^{\frac{i\pi}{1+\alpha}}|$. On the other hand for 
$l> \alpha$ the inequality holds 
\be\label{ineq4}
c_l|1-e^{\frac{i\pi}{l+1}}|^{-1}\gs (c_l l)(l|1-e^{\frac{i\pi}{l+1}}|)^{-1}\gs (c_l l)(\alpha|1-e^{\frac{i\pi}{1+\alpha}}|)^{-1}.
\ee
Using that $c(1)=\sum^{\infty}_{j=0}c_j$ we deduce
\be\label{ineq5}
c_l l\ls \sum^l_{j=1}c_j\ls c(1)-1= \alpha(1-b)
\ee
and thus combining (\ref{ineq4}) and (\ref{ineq5}) we obtain that
$$
c_l|1-e^{\frac{i\pi}{l+1}}|^{-1}\ls (1-b)|1-e^{\frac{i\pi}{1+\alpha}}|^{-1},
$$
which is (\ref{ineq3}). As we have noted the bound
used in (\ref{ineq2}) implies that
$$
|c(e^{i\theta})|\gs (1-b)|1-e^{\frac{\pi i}{1+\alpha}}|^{-1}(\sqrt{1+b(1-b)^{-2}|1-e^{\frac{\pi i}{1+\alpha}}|^2}-1),
$$
which is equivalent to
\be\label{ineq6}
|c(e^{i\theta})|\gs |1-e^{\frac{\pi i}{1+\alpha}}|^{-1}(|(1-b)+b(1-e^{\frac{\pi i}{1+\alpha}})|-(1-b)).
\ee
Plugging (\ref{ineq6}) into (\ref{ineq1}) we derive 
$$
|c(re^{i\theta})|\gs \frac{|(1-b)+b(1-e^{\frac{\pi i}{1+\alpha}})|-(1-b)}{|1-e^{\frac{\pi i}{1+\alpha}}|}-c(r)+c(1).
$$
Finally using (\ref{main0}) we conclude that
$$
\sup_{|z|=r}|\sum^{\infty}_{n=0}(u_n-u_{\infty})z^n|\ls \frac{c(r)-c(1)}{c(1)(r-1)((1-b)D(\alpha)-c(r)+c(1))},
$$ 
where $D(\alpha)=|1-e^{\frac{\pi i}{1+\alpha}}|^{-1}(|1+\frac{b}{1-b}(1-e^{\frac{\pi i}{1+\alpha}})|-1)$
which completes the proof of Theorem \ref{kendal1}.
\end{dwd}

\section{Comparing with the previous bounds} \label{sect3}

\noindent
We recall that our bound on the radius of convergence is of the form $\min_{1\ls \alpha\ls \alpha_0}(1+\frac{D(\alpha)}{\alpha})^{\frac{1}{\kappa(\alpha)}}$.
As it will be shown, this bound is always better than the main bound in \cite{Bax} (Theorem 3.2). Then we turn to study the reason for this improvement. 
Using the limit case with $b,L$ fixed and $R\ra 1$, we check that the minimum of $(1+\frac{D(\alpha)}{\alpha})^{\frac{1}{\kappa(\alpha)}}$ can be attained
in the interval $[1,\alpha_0]$ and that it is data depending problem. On the other hand   
we stress that in the usual setting this minimum of will be attained at $\alpha_0$.
The intuition for this phenomenon is that the smaller is $\pi(C)$ the worse rate of convergence we should expect.
However it is possible to choose $b,R,L$ in a way that this intuition fails and then the minimization procedure is in use. Below we describe how to transform the minimization problem into solving some equation.  
\smallskip

\noindent
Observe that the minimum of the function $(1+\frac{D(\alpha)}{\alpha})^{\frac{1}{\kappa(\alpha)}}$
is attained at the unique point $\alpha$ that satisfies
\be\label{aaa}
\log\frac{N-1}{1-b}=
\log\alpha+\log(1+\frac{D(\alpha)}{\alpha})\frac{D(\alpha)+\alpha}{D(\alpha)-\alpha D'(\alpha)}.
\ee
Obviously to find the minimum on the interval $[1,\alpha_0]$ the solution of (\ref{aaa}) must be compared 
with $1$ and $\alpha_0$. Consequently $r_0=(1+D(1))^{\frac{1}{\kappa(1)}}$ when such $\alpha$ is smaller than $1$ and $r_0=(1+\frac{D(\alpha_0)}{\alpha_0})^{\frac{1}{\kappa(\alpha_0)}}$  when it is bigger than $\alpha_0$, otherwise the solution of (\ref{aaa}) is the worst possible $\alpha$ that minimizes our bound on the radius of convergence. The same discussion concerns maximization of $K_0(r)$. Clearly
the problem reduces to finding the maximum of the function $\alpha (D(\alpha))^{-1}(r^{\kappa(\alpha)}-1)$ which is attained at the unique point $\alpha$ that satisfies the equation
\be\label{a0}
(1+\frac{D'(\alpha)\alpha}{D(\alpha)})(r^{\kappa(\alpha)}-1)=\frac{\log r}{\log R} r^{\kappa(\alpha)}.
\ee
To find the maximum of $\alpha (D(\alpha))^{-1}(r^{\kappa(\alpha)}-1)$
on $[1,\alpha_0]$ we compare the solution of (\ref{a0}) with $1$ and $\alpha_0$. If such $\alpha$ is greater than $\alpha_0$
then 
$$
\alpha_0 (D(\alpha_0))^{-1}(r^{\kappa(\alpha)_0}-1)
$$ 
is the optimal bound on $\max_{1\ls \alpha\ls \alpha_0}\alpha (D(\alpha))^{-1}(r^{\kappa(\alpha)}-1)$.
Similarly if $\alpha\ls 1$ then $(D(1))^{-1}(r^{\kappa(1)}-1)$ is the bound and otherwise
the solution of (\ref{a0}) is the point maximum for $\max_{1\ls \alpha\ls \alpha_0}\alpha (D(\alpha))^{-1}(r^{\kappa(\alpha)}-1)$.
\begin{rema}
It is possible that our Lyapunov function is as good that we predict $R$ as the lower bound on the radius of convergence of $\sum^{\infty}_{n=0}(u_n-u_{\infty})z^n$. Here it is the case when the solution of (\ref{aaa}) 
is smaller than $1$ i.e. when
$$
(1+\frac{D'(1)}{D(1)})(R^{\kappa(1)}-1)\gs R^{\kappa(1)}.
$$
\end{rema}
We turn to show computable bounds on $K_0(r)$ in the case when $\pi(C)$
is unknown. Note that function $D(\alpha)$ is decreasing and therefore $D(\alpha)\gs D(\alpha_0)$. Consequently one can rewrite Corollary \ref{cori2} with $D(\alpha)$ replaced by $D(\alpha_0)$ and in this way obtain new bounds: $K_1(r)\gs K_0(r)$
and $r_1\ls r_0$, where $r_1=\min\{R,\min_{1\ls \alpha\ls \alpha_0}
(1+\frac{D(\alpha_0)}{\alpha})^{\frac{1}{\kappa(\alpha)}}\}$ and
$$
K_1(r)=\max_{1\ls \alpha\ls \alpha_0}\frac{r^{\kappa(\alpha)}-1}{(r-1)(D(\alpha_0)\alpha^{-1}-r^{\kappa(\alpha)+1})}.
$$
Consequently to find $K_1(r)$ it suffices to compute the maximum of $\alpha(r^{\kappa(\alpha)}-1)$ on the interval $[1,\alpha_0]$. The maximum of $\alpha(r^{\kappa(\alpha)}-1)$ is attained at $\alpha$ that satisfies
\be\label{a1}
(r^{\kappa(\alpha)}-1)=\frac{\log r}{\log R} r^{\kappa(\alpha)}.
\ee
There is explicit solution of (\ref{a1}) of the form 
\be\label{a2}
\alpha=\frac{N-1}{1-b}(1-\frac{\log r}{\log R})^{\frac{\log R}{\log r}}.
\ee
Again the solution must be compared with $1$ and $\alpha_0$ which 
finally provides the direct form of $K_1(r)$. 
\begin{coro}\label{cori3}  
Suppose that $b_1\gs b$ and $b(R)\ls L$. 
\begin{enumerate}
\item
If $1\gs \frac{N-1}{1-b}(1-\frac{\log r}{\log R})^{\frac{\log R}{\log r}}$, then
$$
\sup_{|z|=r}|\sum^{\infty}_{n=0}(u_n-u_{\infty})z^n|\ls K_1(r)= (r-1)^{-1}([\frac{D(\alpha_0)}{(r^{\kappa(1)-1})}-1]_{+})^{-1}. 
$$
\item
If $1\ls \frac{N-1}{1-b}(1-\frac{\log r}{\log R})^{\frac{\log R}{\log r}}\ls \alpha_0$, then
\begin{align*}
&\sup_{|z|=r}|\sum^{\infty}_{n=0}(u_n-u_{\infty})z^n|\ls K_1(r)=\\ &=(r-1)^{-1}([\frac{(1-b)D(\alpha_0)}{N-1}\frac{\log R}{\log r}(1-\frac{\log r}{\log R})^{-\frac{\log R}{\log r}+1}-1]_{+})^{-1}.
\end{align*}
\item
If $\alpha_0\ls \frac{N-1}{1-b}(1-\frac{\log r}{\log R})^{\frac{\log R}{\log r}}$, then
$$
\sup_{|z|=r}|\sum^{\infty}_{n=0}(u_n-u_{\infty})z^n|\ls K_1(r)= 
(r-1)^{-1}([\frac{D(\alpha_0)}{\alpha_0(r^{\kappa(\alpha_0)-1})}-1]_{+})^{-1}.
$$
\end{enumerate}
\end{coro}
Corollary \ref{cori3} implies the following expression for $r_1$. Let $x_{\alpha}=r$, $\alpha\gs 0$ be the unique solution of 
\be\label{origami}
\alpha=\frac{N-1}{1-b}(1-\frac{\log r}{\log R})^{\frac{\log R}{\log r}}
\ee 
if $\frac{N-1}{(1-b)\alpha}\gs e$ and $x_{\alpha}=1$ otherwise. 
\begin{coro}\label{cori5}
Suppose that $b_1\gs b$ and $b(R)\ls L$. Let $\bar{r}$ be the unique solution of 
$$
\frac{(1-b)D(\alpha_0)}{N-1}=\frac{\log r}{\log R}(1-\frac{\log r}{\log R})^{\frac{\log R}{\log r}-1},
$$
if $\bar{r}\ls x_1$ then $r_1=(1+D(\alpha_0))^{\frac{1}{\kappa(1)}}$, if $x_1\ls \bar{r}\ls x_{\alpha_0}$
then $r_1=\bar{r}$ and if $\bar{r}\gs x_{\alpha_0}$ then $r_1=(1+\frac{D(\alpha_0)}{\alpha_0})^{\frac{1}{\kappa(\alpha_0)}}$.
\end{coro}
Clearly $r_1\gs r_0$, we turn to show that $r_1$ is a better bound than 
in Theorem 3.2 in \cite{Bax}. For this reason denote by $r_2$ the unique solution of
\be\label{nif2}
\frac{r-1}{r}\frac{1}{\log^2(R/r)}=\frac{b}{2N}.
\ee
Our aim is to show that $r_2\ls r_1$. 
First observe that by the definition 
$$
r^{\kappa(\alpha)}_1-1=\frac{D(\alpha_0)}{\alpha},
$$
for some $\alpha\in [1,\alpha_0]$.
Again by the definition $R^{\kappa(\alpha)}=(N-1)/((1-b)\alpha)$, which implies that 
\be\label{pop}
\kappa(\alpha)r^{\kappa(\alpha)}_1\frac{r_1-1}{r_1}\gs r^{\kappa(\alpha)}-1\gs \frac{(1-b) R^{\kappa(\alpha)} D(\alpha_0)}{N-1}.
\ee
By the following inequality 
$$
D(\alpha_0)=\frac{\sqrt{(1-b)^2+4b\sin^2(\frac{\pi}{2(1+\alpha_0)})}-(1-b)}{2(1-b)\sin(\frac{\pi}{2(1+\alpha_0)})}\gs \frac{b}{(1-b)(1+\alpha_0)},
$$
we obtain in (\ref{pop}) 
\be\label{pop1}
\kappa(\alpha)\frac{r^{\kappa(\alpha)}_1}{R^{\kappa(\alpha)}}\frac{r_1-1}{r_1}\gs \frac{b}{(1+\alpha_0)(N-1)}.
\ee
It suffices to note that $(1+\alpha_0)\ls 2\kappa(\alpha_0)\ls 2\kappa(\alpha)$, which is the consequence of $\kappa(\alpha_0)\ls \kappa(\alpha)$ and the fact that 
$$
R^{\kappa(\alpha_0)}=R\frac{R^{\alpha_0}-1}{R-1},
$$
which can be used to show that for a given $R$, the function $\kappa(\alpha_0)/(1+\alpha_0)$ is increasing with $\alpha_0$. Thus since $\kappa(\alpha_0)/(1+\alpha_0)= 1/2$ for $\alpha_0=1$
we deduce that $(1+\alpha_0)\ls 2\kappa(\alpha_0)$. 
Plugging the estimate $2\kappa(\alpha)\gs (1+\alpha_0)$
into (\ref{pop1}) we derive
$$
\kappa(\alpha)^2\frac{r^{\kappa(\alpha)}_1}{R^{\kappa(\alpha)}}\frac{r_1-1}{r_1}\gs \frac{b}{2(N-1)}.
$$
It remains to check that $\kappa(\alpha)=2/\log(R/r_1)$
is the point maximum of $\kappa(\alpha)^2 (r_1/R)^{\kappa(\alpha)}$, which follows that 
$$
\frac{r_1}{r_1-1}\frac{1}{\log^2(R/r_1)}\gs \frac{b e^2}{8(N-1)}.
$$
This shows that $r_1\gs r_2$ and in fact $r_2$ can be treated as the lower bond in the worst possible case of our result.
We stress that using $\alpha_0$ instead of the minimization over all $\alpha_0$ usually gives a major numerical improvement. 
\smallskip

\noindent
To provide a rough argument for exploiting the parameter $\alpha_0$ we consider the simplest renewal model where 
there are only two possible states $1$ and $\alpha_0$ (for simplicity assume that $\alpha_0\in\N$). Then  
the optimal rate of convergence is closely related to the specific solution of 
$\frac{bz+(1-b)z^{\alpha_0}-1}{z-1}=0$, namely it is the inverse of the smallest absolute value of the solution of the equation. 
Denoting the root by $z_0$ one can show that
\be\label{odc4}
|z_{\alpha_0}|= 1+\frac{2b\pi^2}{(1-b)^2\alpha^3_0}+o(\alpha^{-3}_0),
\ee
(see discussion after Theorem 3.2 in \cite{Bax}) and $\alpha_0$ is exactly our parameter. 
Therefore whenever the estimate $(1+\frac{D(\alpha_0)}{\alpha_0})^{\frac{1}{\kappa(\alpha_0)}}$ is applied one cannot improve it 
up to numerical constant. We turn to study this phenomenon in the limit case where $b,L$ are fixed and $R\ra 1$.
\begin{coro}\label{cori4}
Suppose that $R\ra 1$ and $b_1\gs b$, $b(R)\ls L$.
\begin{enumerate}
\item If $(\frac{L-1}{1-b})/\log \frac{L-b}{1-b}\gs e^{1/2}$, then
$$
r_0(R) =1+ \frac{b\pi (R-1)^3}{2(1-b)^2}\log^{-2}(\frac{L-b}{1-b})\log^{-1}((L-1)/\log \frac{L-b}{1-b})+o((R-1)^3),
$$
\item If $(\frac{L-1}{1-b})/\log \frac{L-b}{1-b}\ls e^{1/2}$, then
$$
r_0(R) =1+ \frac{b e\pi (R-1)^3}{(L-1)^2}+o((R-1)^3), 
$$
\end{enumerate}
\end{coro}
\begin{dwd} 
Observe that $\lim_{\alpha\ra\infty}\alpha D(\alpha)=\frac{b\pi}{2(1-b)^2}$, thus
we can treat $\pi b(2(1-b)^2\alpha)^{-1}$ as the right approximation of $D(\alpha)$ when $\alpha$ 
tends to infinity. As we have stated in Corollary \ref{cori2} to find
\be\label{redek}
r_0(R)=\inf_{1\ls \alpha\ls \alpha_0(R)}(1+\frac{D(\alpha)}{\alpha})^{\frac{1}{\kappa(\alpha)}}
\ee
one should solve the equation (\ref{aaa}), i.e. find $\alpha(R)$ that satisfies
\be\label{redek2}
\log\frac{N(R)-1}{1-b}=
\log\alpha+\log(1+\frac{D(\alpha)}{\alpha})\frac{D(\alpha)+\alpha}{D(\alpha)-\alpha D'(\alpha)},
\ee
where $N(R)=(L-1)/(R-1)$, and compare the outcome with $1$ and $\alpha_0(R)$. In particular we deduce from (\ref{redek2}) that $\alpha(R)$ necessarily tends to infinity when $R\ra 1$,
hence using
$$
\lim_{\alpha\ra \infty}(1+\frac{\alpha}{D(\alpha)})\log(1+\frac{D(\alpha)}{\alpha})=1\;\;\mbox{and}\;\;
\lim_{\alpha\ra \infty}(1-\frac{\alpha D'(\alpha)}{D(\alpha)})=2,
$$
we obtain that
$$
\log \alpha(R)=-\frac{1}{2}+\log \frac{N(R)-1}{1-b}+o(1). 
$$
The solution must be compared with $\alpha_0(R)$ therefore if 
$$
\lim_{R\ra \infty} \frac{N(R)-1}{\alpha_0(R)R-1}=\frac{L-1}{1-b}\log^{-1}(\frac{L-b}{1-b})< e^{\frac{1}{2}} 
$$
we have to use $\alpha(R)$ (at least for small enough $R$) when minimize $(1+\frac{D(\alpha)}{\alpha})^{\frac{1}{\kappa(\alpha)}}$ over $[1,\alpha_0(R)]$,
otherwise $\alpha_0(R)$ is the point minimum.
In the first setting we have
$$
\alpha(R)=e^{-\frac{1}{2}}\frac{L-1}{(1-b)(R-1)}+o(1),\;\;\mbox{and}\;\;\kappa(\alpha(R))=
\frac{1}{2(R-1)}+o(1),
$$
thus using (\ref{redek}) we obtain that
\begin{align*}
& r_0(R)=(1+\frac{D(\alpha(R))}{\alpha(R)})^{\frac{1}{\kappa(\alpha(R))}}=1+\frac{D(\alpha(R))}{\alpha(R)\kappa(\alpha(R))}+o(R-1)=1+\\
& +\frac{\pi b}{2(1-b)^2\alpha^2(R)\kappa(\alpha(R))}+o((R-1))=1+\frac{\pi e b (R-1)^3 }{(L-1)^2}+o((R-1)^3).
\end{align*}
In the same way if $\frac{L-1}{1-b}\log^{-1}(\frac{L-b}{1-b})\gs e^{1/2}$, then
$$
\alpha_0(R)=\frac{\log(\frac{L-b}{1-b})}{R-1}+o(1),\;\;\kappa(\alpha_0(R))=
\frac{\frac{L-1}{1-b}}{(R-1)\log(\frac{L-b}{1-b})}+o(1), 
$$
and hence
\begin{align*}
& r_0(R)=(1+\frac{D(\alpha_0(R))}{\alpha_0(R)})^{\frac{1}{\kappa(\alpha_0(R))}}=1+\frac{D(\alpha_0(R))}{\alpha_0(R)\kappa(\alpha_0(R))}+o(R-1)=\\
& =1+\frac{\pi b}{2(1-b)^2\alpha^2_0(R)\kappa(\alpha_0(R))}+o((R-1))=\\ 
& =1+\frac{\pi b (R-1)^3 }{2(1-b)^2}\log^{-2}(\frac{L-b}{1-b})\log^{-1}((L-1)/\log(\frac{L-b}{1-b}))+o((R-1)^3).
\end{align*}
It completes the proof of the corollary.  
\end{dwd}

\noindent
In particular Corollary \ref{cori4} shows that whenever  
$\frac{L-1}{1-b}\log^{-1}(\frac{L-b}{1-b})\gs e^{\frac{1}{2}}$ the
following inequality holds
$$
r_0(R)=1+\frac{\pi b}{2(1-b)^2 \alpha_0^2(R) \kappa(\alpha_0(R))}+o(\alpha_0(R)^{-3}),
$$
which when compared with (\ref{odc4}) proves that our result cannot be improved up to a numerical constant (we recall that $(1+\alpha_0)/2\ls \kappa(\alpha_0)\ls \alpha_0$). 
On the other hand Corollary \ref{cori4} makes it possible to compare  
our result with Theorem 3.2 in \cite{Bax}. 
The following estimate holds for $r_2$ (see Section 3 in \cite{Bax}) 
$$
r_2(R)=1+\frac{e^2 b (R-1)^3}{8(L-1)}+o((R-1)^3).
$$
Therefore if $L-1$ much larger than $1-b$ our answer is better on $(L-1)/(1-b)^2$ and if $L-1$ is close to $1-b$ then 
on $L-1$. 
\smallskip

\noindent
Nevertheless there are indeed two cases depending on $b,R,L$ either $L$ is far from $1$ with respect to $b,L$ and then 
the minimum of $(1+\frac{D(\alpha)}{\alpha})^{\frac{1}{\kappa(\alpha)}}$ is attained on $\alpha_0(R)$.
However it may happen that $L$ is close to $1$ (again with respect to $b$ and $L$) and then
we have to use the minimization inside $[1,\alpha_0(R)]$ even for $R\ra 1$. It explains that the minimization of our bound on $r_0$
for $\alpha\in [1,\alpha_0]$ is important in our discussion.

\section{Applications to the geometric ergodicity}\label{sect4}

In this section we follow our discussion from the introduction.
We start from the atomic case, note that Proposition \ref{pro2} is the classical first entrance last exit decomposition rule considered in two cases in order to obtain better bounds on $M_W(r)$. To show this
it suffices to carefully read the proof of Proposition 4.2 in \cite{Bax}.
In Theorem \ref{thm1} we have stated the best form of such an estimate in case of
$W\equiv 1$. To obtain the result we need a small modification of Proposition 4.2 \cite{Bax}, we give the 
argument for the sake of completeness.
\begin{dwd}[of Proposition \ref{pro3}]
To show the first inequality it suffices to observe that $r^{-1}H_1(r,x)$ attains
its maximum on the interval $[1,\lambda^{-1}]$ at $\lambda^{-1}$. Using Proposition \ref{pro1}
we obtain that 
$$
r^{-1}H_1(r,x)\ls \lambda H_1(\lambda^{-1},x)=\frac{G(\lambda^{-1},x)-1}{\lambda^{-1}-1}\ls \frac{V(x)-1}{\lambda^{-1}-1}.
$$
Consequently $H_1(r,x)\ls \frac{r\lambda(V(x)-1)}{1-\lambda}$, $x\not\in C$  and 
in the same way we show that $H_1(r,x)\ls \frac{r(K-\lambda)}{1-\lambda}$ if $x\in C$.
The second inequality can be derived in the similar way, first we note that 
$r^{-1}(r-1)^{-1}(H_1(r,a)-rH_1(1,a))$ is increasing and then we use the bound
$$
\lambda \frac{H_1(\lambda^{-1},a)-\lambda^{-1}H_1(1,a)}{\lambda^{-1}-1}\ls 
\frac{\frac{K-\lambda}{1-\lambda}-1}{1-\lambda}=\frac{K-1}{(1-\lambda)^2}.
$$
\end{dwd}

\noindent
Theorems \ref{thm1} and \ref{thm2} are clear consequence of Proposition \ref{pro2} and the bounds
we have proved, which completes the atomic part of the discussion. 
\smallskip

\noindent
In the general case we do not require $C$ to be an atom. 
However, there is a simple trick (cf. Meyn -Tweedie \cite{MT1}, Numellin \cite{Num}) which reduces this case to the atomic one. Consider the split chain $(X_n,Y_n)_{n\gs 0}$ defined on state space $\bar{\ccS}=\ccS\times\{0,1\}$ with the $\sigma$-field $\bar{\ccB}$ generated by $\ccB\times \{0\}$ and $\ccB\times\{1\}$. We define transition probabilities as follows:
\begin{align*}
& \P(Y_n=1| \ccF^{X}_{n}, \ccF^{Y}_{n-1})=\bar{b}1_{C}(X_n),\\
& \P(X_{n+1}\in A|\ccF^X_{n},\ccF^{Y}_{n})=\left\{\begin{array}{ll}
\nu(A), & \mbox{if}\;Y_n=1,\\
 \frac{\P(X_n,A)-\bar{b}1_C(X_n)\nu(A)}{1-\bar{b}1_C(X_n)}, & \mbox{if} \;Y_n=0.
\end{array}\right.
\end{align*}
where $\ccF^X_n=\sigma(X_k:\;0\ls k\ls n)$, $\ccF^Y_n=\sigma(Y_k:\;0\ls k\ls n)$. Thus the chain evolves in a way that whenever
$X_n$ is in $C$ we pick $Y_{n}=1$ with probability $\bar{b}$. Then if $Y_n=1$ we chose $X_{n+1}$ from $\nu$ distribution whereas if
$Y_n=0$ then we just apply normalized probability measure version of $\P(X_n,\cdot)-\bar{b}1_C\nu$. The split chain is designed so that it has an atom $\ccS\times\{1\}$
and so that its first component $(X_n)_{n\gs 0}$ is a copy of the original Markov chain. Therefore we can apply the approach from the previous section to
the split chain $(X_n,Y_n)$ and the stopping time
$$
T=\min\{n\gs 1:\;Y_n=1\}.
$$
Let $\P_{x,i}$, $\E_{x,i}$ denote the probability and the expectation for the split chain started with $X_0=x$ and $Y_0=i$. Observe that for a fixed point $a\in C$
we have $\P_{x,1}=\P_{a,1}$ and $\E_{x,1}=\E_{a,1}$ for all $x\in C$. 
Following the method used in the atomic case we define the renewal sequence $\bar{u}_n=\P_{a,1}(Y_n=1)$ and the corresponding increment sequence $\bar{b}_n=\P_{a,1}(T=n)$ for $n\gs 1$. Clearly $\bar{u}_n=\P_{a,1}(X_n\in C, Y_n=1)=\bar{b}\P_{\nu}(X_{n-1}\in C)$ for $n\gs 1$, so
\be\label{ark1}
\bar{b}_1=\bar{b}\nu(C)\gs b,\;\;\mbox{and}\;\;\bar{u}_{\infty}=\bar{b}\pi(C).
\ee
We define
$$
\bar{G}(r,x,i):=\E_{x,i}(r^T),\;\;\bar{H}_W(r,x,i):=\E_{x,i}(\sum^{T}_{n=1}r^n W(X_n)),
$$
for all $x\in \ccS$, $i=0,1$ and all $r>0$ for which the right hand sides are well defined. 
We also need the following expectation 
$$
\E_x:=(1-\bar{b}1_C(x))\E_{x,0}+\bar{b}1_C(x)\E_{x,1},
$$
which agrees with the usual $\E_x$ on $\ccF^X$.
There exists a unique stationary measure $\bar{\pi}$ say, on $(\bar{\ccS},\bar{\ccB})$, so that $\int \bar{g}\bar{\pi}=\int g d\pi$ (where $g(x)=\bar{g}(x,0)=\bar{g}(x,1)$ for all $x\in \ccS$). In particular we have that $\bar{\pi}(\ccS\times\{1\})=\bar{b}\pi(C)$.
The first entrance last exist decomposition leads to the following result (cf. 
Proposition \ref{pro2} and Proposition 4.3 in \cite{Bax})
\begin{pros}\label{pro5}
For all $a\in \ccS\times \{1\}$
\begin{eqnarray}
&&\sup_{|z|= r}|\sum^{\infty}_{n=1}(P^n \bar{g}(a)-\int g d\pi)z^n|\ls \bar{H}_W(r,a,1)\sup_{|z|= r}|\sum^{\infty}_{n=0}(\bar{u}_n-\bar{u}_{\infty})z^n|+\nonumber\\
\label{bam1}&&+\bar{b}\pi(C)\frac{\bar{H}_W(r,a,1)-r\bar{H}_W(1,a,1)}{r-1},
\end{eqnarray}
and for all $x \in \ccS\times \{0\}$ 
\begin{eqnarray}
&&\sup_{|z|= r}|\sum^{\infty}_{n=1}(P^n \bar{g}(x)-\int g d\pi)z^n|\ls\nonumber\\
&&\ls \bar{H}_W(r,x,0)+
\bar{G}(r,x,0)\bar{H}_W(r,a,1)\sup_{|z|= r}|\sum^{\infty}_{n=0}(\bar{u}_n-\bar{u}_{\infty})z^n|+\nonumber\\
\label{bam2}&&+\bar{b}\pi(C)\frac{\bar{H}_W(r,a,1)-r\bar{H}_W(1,a,1)}{r-1}\bar{G}(r,x,0)+\nonumber\\
&&+\bar{b}\pi(C)\bar{H}_W(1,a,1)\frac{r(\bar{G}(r,x,0)-1)}{r-1}. 
\end{eqnarray}
\end{pros}
Therefore one has to estimate all quantities in the above result. We move the exact
calculations to the appendinx as well as the examples with computations of the convergence rate for some toy examples.

\section*{Appendix A}\label{appendA}

\renewcommand{\thethm}{A.\arabic{thm}}
\setcounter{thm}{0}

\renewcommand{\thesubsection}{A.\arabic{subsection}}
\setcounter{subsection}{0}

\renewcommand{\theequation}{A.\arabic{equation}}
\setcounter{equation}{0}

\subsection {Global bounds}

Our method described in Corollary \ref{cori6} implies that 
$$
\sup_{|z|= r}|\sum^{\infty}_{n=0}(\bar{u}_n-\bar{u}_{\infty})z^n|\ls K_0(r)\;\;\mbox{for}\; 1\ls r\ls r_0.
$$
The first step is to replace 
the stopping time $T$ by $\tau=\tau_C$. For this reason we define
$$
G(r,x,i)=\E_{x,i}r^{\tau},\;\;H_W(r,x,i)=\E_{x,i}(\sum^{\tau}_{n=1}r^n W(X_n)).
$$
Let also $G(r)=\sup_{x\in C}\E_{x,0}r^{\tau}$, $H_W(r)=\sup_{x\in C}\E_{x,0} \sum^{\tau}_{n=1}r^{n}W(X_n)$. In the Lemma A.1 in \cite{Bax} there are proved following inequalities:
\begin{pros}\label{pro6}
For $r\ls \lambda^{-1}$ and $(1-\bar{b})G(r)<1$ the inequalities hold:
\be\label{onty1}
\bar{G}(r,x,i)\ls \frac{\bar{b}G(r,x,i)}{1-(1-\bar{b})G(r)}
\ee
and
\be\label{onty2}
\bar{H}_W(r,x,i)\ls H_W(r,x,i)+\frac{(1-\bar{b})H_W(r)G(r,x,i)}{1-(1-\bar{b})G(r)}
\ee
\end{pros}
In the introduction we have explained that the crucial for our approach is to establish (\ref{ing1}).
We have all necessary tools to get the result.
\begin{pros}\label{pro7}
For all $a\in C$ and $1\ls r\ls \min\{\lambda^{-1},(1-\bar{b})^{-\frac{1}{1+\alpha_1}}\}$
\be\label{onty3}
\bar{G}(r,a,1)\ls \max\{\frac{\bar{b}r}{1-(1-\bar{b})r^{1+\alpha_1}},\frac{br+(\bar{b}-b)r^{\alpha_2}}{1-(1-\bar{b})r}\},
\ee
where $\alpha_1=\log(\frac{(K-\bar{b})}{(1-\bar{b})})/\log \lambda^{-1}$, $\alpha_2=\log(\frac{(K-1+\bar{b}-b)}{\bar{b}-b})/\log \lambda^{-1}$.
Moreover if $1+b\gs 2\bar{b}$, then
\be\label{onty35}
\bar{G}(r,a,1)\ls \frac{\bar{b}r}{1-(1-\bar{b})r^{1+\alpha_1}}.
\ee
For all $x\in \ccS$, $1\ls r\ls \min\{\lambda^{-1},(1-\bar{b})^{-\frac{1}{1+\alpha_1}}\}$
\be\label{onty4}
\bar{b}1_{C}(x)+(1-\bar{b}1_{C}(x))\bar{G}(r,x,0)\ls 
\frac{\bar{b}V(x)}{1-(1-\bar{b})r^{1+\alpha_1}}.
\ee
\end{pros}
\begin{dwd}
The split chain construction implies that for any $a\in C$
\be\label{superlaska}
(1-\bar{b})\sup_{x\in C}G(r,x,0)+\bar{b}G(r,a,1)=\sup_{x\in C}G(r,x)=G(r).
\ee
Moreover due to $\bar{b}\nu(C)\gs b$ we have that $\bar{b}G(r,a,1)=\bar{b}\sum^{\infty}_{k=1}\P_{\nu}(\sigma=k-1)r^k$, where $\sigma=\inf\{n\gs 0:\;X_n\in C\}$
has its first coefficient greater or equal $b$. Therefore by our usual argument with the H\"{o}lder inequality we deduce that
$$
(1-\bar{b})\sup_{x\in C}G(r,x,0)\ls (1-\bar{b})rv^{\frac{\log r}{\log \lambda^{-1}}},\;\;\mbox{and},\;\;\bar{b}G(r,a,1)\ls br+(\bar{b}-b)r u^{\frac{\log r}{\log \lambda^{-1}}},
$$
where $u=G(\lambda^{-1},a,1),v=\sup_{x\in C}G(\lambda^{-1},x,0)$ verify
\be\label{baff}
\lambda^{-1}(b+(\bar{b}-b)u+(1-\bar{b})v)=\sup_{x\in C}G(\lambda^{-1},x)\ls K\lambda^{-1},\;\;u,v\gs 1.
\ee
Observe that by (\ref{onty1}) we have the following bound 
\be\label{baff1}
\bar{G}(r,a,1)\ls \frac{\bar{b}G(r,a,1)}{1-(1-\bar{b})G(r)}\ls F(u,v)=\frac{br+(\bar{b}-b)ru^{\frac{\log r}{\log \lambda^{-1}}})}{1-(1-\bar{b})rv^{\frac{\log r}{\log \lambda^{-1}}}}.
\ee
One can check that the bounding function $F(u,v)$ is convex for all $(u,v)$
that satisfy (\ref{baff}) and hence it takes its maximum on the boundaries of the set given by (\ref{baff}). Consequently due to (\ref{baff1})
we obtain that
\be\label{baff2}
\bar{G}(r,a,1)\ls \max\{\frac{\bar{b}r}{1-(1-\bar{b})r^{1+\alpha_1}},\frac{br+(\bar{b}-b)r^{1+\alpha_2}}{1-(1-\bar{b})r}\}.
\ee
It is easy to check that whenever $1+b\gs 2\bar{b}$ then $(1-\bar{b})\alpha_1\gs (\bar{b}-b)\alpha_2$ and
the maximum in (\ref{baff2}) can be replaced by the first quantity for any $r\gs 1$.
Otherwise if $1+b<2\bar{b}$ then $(1-\bar{b})\alpha_1<(\bar{b}-b)\alpha_2$ and therefore for small enough $r$ the maximum in (\ref{baff2}) 
must be attained at the second expression.
\smallskip

\noindent
We turn to show the second assertion. Observe that by Proposition \ref{pro1} we have  
$G(r,x,0)=G(r,x)\ls V(x)$ for all $x\not\in C$. Consequently (\ref{onty1}) yields
\be\label{ag4}
\bar{G}(r,x,0)\ls \frac{\bar{b}V(x)}{1-(1-\bar{b})G(r)},
\ee
for all $x\not\in C$. Since obviously $G(r)\ls r^{1+\alpha_1}$ we deduce that
$$
\bar{G}(r,x,0)\ls \frac{\bar{b}V(x)}{1-(1-\bar{b})r^{\alpha_1}}
$$
On the other hand by (\ref{onty1})
$$
\bar{G}(r,x,0)\ls \frac{\bar{b}r^{\alpha_1}}{1-(1-\bar{b})r^{1+\alpha_1}}
$$
for all $x\in C$ and therefore
\be\label{ag5}
\bar{b}+(1-\bar{b})\bar{G}(r,x,0)\ls \frac{\bar{b}}{1-(1-\bar{b})r^{1+\alpha_1}}
\ee
for all $x\in C$. Since $V\gs 1$, inequalities (\ref{ag4}) and (\ref{ag5}) imply (\ref{onty4}).
\end{dwd}
The next step is to obtain the estimate (\ref{antyma3}). 
\begin{coro}\label{jasmine2}
The following inequality holds
\be\label{kotek1}
\bar{b}^{-1}\ls \frac{\bar{H}_1(1,a,1)-1}{1-b}\ls \bar{b}^{-1}\max\{\frac{1-\bar{b}}{1-b}(1+\alpha_1),
\frac{1-\bar{b}}{1-b}+\frac{\bar{b}-b}{1-b}\alpha_2\}=\bar{b}^{-1}\bar{\alpha}_0.
\ee
\end{coro}
\begin{dwd}
As for the first assertion we simply apply (\ref{onty3}) to
bound $\bar{H}_1(r,a,1)=\frac{r\bar{G}(r,a,1)-1}{r-1}$ and then
tend with $r$ to $1$. To prove the second assertion let
$S=\max\{k\gs 1:\tau_k\ls T\}$, where $\tau_k$, $k\gs 0$ are subsequent visits to $C$ by $(X_n)_{n\gs 0}$, in particular $\tau_0=0$. Observe that 
$$
\bar{H}_1(1,a,1)=\E_{a,1}(\sum^{\infty}_{k=0}1_{S\gs k}(\tau_k-\tau_{k-1})).
$$
Therefore by the construction
$$
\bar{H}_1(1,a,1)\gs \E_{\nu}(1+\sigma)+\E_{a,1}(S-1),
$$
where we recall that $\sigma=\min\{n\gs 0:\;X_n\in C\}$. 
Since $S$ has the geometric distribution with the probability of success
$\bar{b}$ we obtain that
$$
\bar{H}_1(1,a,1)\gs \bar{b}^{-1}+\E_{\nu}\sigma.
$$
It remains to notice that $\E_{\nu}\sigma\gs 1-\nu(C)$, therefore
if $\bar{b}\nu(C)=b$ then 
$$
\bar{H}_1(1,a,1)\gs \bar{b}^{-1}+1-\frac{b}{\bar{b}},
$$
which completes the proof.
\end{dwd}
Now we state an improvement of the result  mentioned in the proof of Proposition 4.4 in \cite{Bax}. 
\begin{pros}\label{pro61}
For $r\ls \lambda^{-1}$ and $(1-\bar{b})G(r)<1$ we have that
\be\label{onty11}
\bar{H}_W(r,a,1)\ls \frac{1}{\bar{b}}\sup_{x\in C}H_W(r,x)+\frac{1-\bar{b}}{\bar{b}}\frac{H_W(r)\sup_{x\in C}(G(r,x)-1)}{1-(1-\bar{b})G(r)} 
\ee
and
\begin{eqnarray}
&&\bar{H}_W(r,a,1)-r\bar{H}_W(1,a,1)\ls
\frac{1}{\bar{b}}\sup_{x\in C}(H_W(r,x,0)-rH_W(1,x,0))+\nonumber\\
\label{onty21}&& +\frac{1-\bar{b}}{\bar{b}}H_W(r)(\bar{G}(r,a,1)-1).
\end{eqnarray}
\end{pros}
\begin{dwd}
To prove the first assertion note that (\ref{superlaska})
can be rewritten as
$$
\frac{\bar{b}G(r,a,1)}{1-(1-\bar{b})G(r)}\ls 1+\frac{\sup_{x\in C}(G(r,x)-1)}{1-(1-\bar{b})G(r)}.
$$
Combining the above inequality with (\ref{onty2}) we derive
$$
\bar{H}_W(r,a,1)\ls H_W(r,a,1)+\frac{1-\bar{b}}{\bar{b}}H_W(r)+\frac{(1-\bar{b})H_W(r)\sup_{x\in C}(G(r,x)-1)}{\bar{b}(1-(1-\bar{b})G(r))}.
$$
Since the definition of $H_W(r,x,1)$ implies that 
$$
\bar{b}H_W(r,a,1)+(1-\bar{b})H_W(r)\ls \sup_{x\in C}H_W(r,x)
$$
we obtain (\ref{onty11}). To show the second assertion we use $S=\max\{k\gs 1:\;\tau_k\ls T\}$ defined in the proof of Corollary \ref{jasmine2}. 
The following inequality holds
\begin{eqnarray}
&&\bar{H}_W(r,a,1)-r\bar{H}_W(1,a,1)\ls H_W(r,a,1)-rH_W(1,a,1)+\nonumber\\
\label{pifpaf}&&+\sum^{\infty}_{k=2}\E_{a,1} [1_{k\ls S}\sup_{x\in C}(r^{\tau_{k-1}}H_W(r,x,0)-rH_W(1,x,0))]. 
\end{eqnarray}
As we have shown in Corollary \ref{jasmine2}, $\E_{a,1}(S-1)= (1-\bar{b})/\bar{b}$ we deduce 
that
\begin{align*}
& \sum^{\infty}_{k=2}(\E_{a,1} 1_{k\ls N})\sup_{x\in C}(H_W(r,x,0)-rH_W(1,x,0))=\\
&=\frac{1-\bar{b}}{\bar{b}}\sup_{x\in C}(H_W(r,x,0)-rH_W(1,x,0)),
\end{align*}
which together with (\ref{pifpaf}) provides
\begin{eqnarray}
&& \bar{H}_W(r,a,1)-r\bar{H}_W(1,a,1)\ls H_W(r,a,1)-rH_W(1,a,1)+\nonumber\\
&&+\frac{1-\bar{b}}{\bar{b}}\sup_{x\in C}(H_W(r,x,0)-rH_W(1,x,0))+\\
\label{ddt}&&+\sum^{\infty}_{k=2}[\E_{a,1}1_{k\ls N}(r^{\tau_{k-1}}-1)]\sup_{x\in C} H_W(r,x,0).
\end{eqnarray}
As usual we observe that
\begin{eqnarray}
&& H_W(r,a,1)-rH_W(1,a,1)+\nonumber\\
&&+\frac{1-\bar{b}}{\bar{b}}\sup_{x\in C}(H_W(r,x,0)-rH_W(1,x,0))\ls \nonumber\\
\label{ftj1} && \ls \frac{1}{\bar{b}}\sup_{x\in C}(H_W(r,x)-rH_W(1,x)).
\end{eqnarray}
Moreover since $Y_{\tau_k}$ is independent of $\tau_{k-1}$ we have  $\E_{a,1} r^{\tau_{k-1}}1_{k\ls S}=(1-\bar{b})\E_{a,1} r^{\tau_{k-1}}1_{k-1\ls S}$ which implies that
$$
\sum^{\infty}_{k=2}\E_{a,1} r^{\tau_{k-1}}1_{S=k-1}=\bar{b}\sum^{\infty}_{k=2}\E_{a,1} r^{\tau_{k-1}}1_{k-1\ls S}=\frac{b}{1-\bar{b}}\sum^{\infty}_{k=2}\E r^{\tau_{k-1}}1_{k\ls S},
$$
and thus
$$
\sum^{\infty}_{k=2}[\E_{a,1}1_{k\ls S}(r^{\tau_{k-1}}-1)]= \frac{1-\bar{b}}{\bar{b}}(\E_{a,1}(r^T-1))=\frac{1-\bar{b}}{\bar{b}}(\bar{G}(r,a,1)-1).
$$
Consequently
\be\label{sgh}
\sum^{\infty}_{k=2}[\E_{a,1}1_{k\ls N}(r^{\tau_{k-1}}-1)]\sup_{x\in C} H_W(r,x,0)=
\frac{1-\bar{b}}{\bar{b}}H_W(r)(\bar{G}(r,a,1)-1).
\ee
Combining (\ref{ddt}), (\ref{ftj1}) and (\ref{sgh})
we conclude that
\begin{align*}
& \bar{H}_W(r,a,1)-r\bar{H}_W(1,a,1)\ls
\frac{1}{\bar{b}}\sup_{x\in C}(H_W(r,x,0)-rH_W(1,x,0))+\\
&+\frac{1-\bar{b}}{\bar{b}}H_W(r)(\bar{G}(r,a,1)-1).
\end{align*}
It completes the proof of (\ref{onty21}).
\end{dwd}

\subsection{Case of $W\equiv 1$}

In the case of $W\equiv 1$, the above result gives 
an improvement in $
\bar{H}_1(r,a,1)-r\bar{H}_1(1,a,1)$  estimation, which as we have mentioned in the introduction,
can be used in the part of the proof where $\sup_{|z|= r}|\sum^{\infty}_{n=0}(\bar{u}_n-\bar{u}_{\infty})|$
is considered.
\begin{pros}\label{pro8}
The following inequalities hold
\be\label{onty5}
(1-\bar{b}1_C(x))\bar{H}_1(r,x,0)\ls \frac{r\lambda(V(x)-1)}{1-\lambda}+
\frac{(1-\bar{b})(r^{1+\alpha_1}-1)rV(x)}{(r-1)(1-(1-\bar{b})r^{1+\alpha_1})},
\ee
for all $x\in \ccS$, $1\ls r\ls \min\{\lambda^{-1},(1-\bar{b})^{-\frac{1}{1+\alpha_1}}\}$,
\be\label{onty6}
\bar{H}_1(r,a,1)\ls   \frac{1}{1-(1-\bar{b})r^{1+\alpha_1}}\frac{r(K-\lambda)}{1-\lambda}
\ee
for all $a\in C$, $1\ls r\ls \min\{\lambda^{-1},(1-\bar{b})^{-\frac{1}{1+\alpha_1}}\}$ and  
\begin{eqnarray}
&&\frac{\bar{H}_1(r,a,1)-r\bar{H}_1(r,a,1)}{r-1}\ls \frac{1}{\bar{b}}\frac{r\lambda(K-1)}{(1-\lambda)^2}+\nonumber\\
&&+\label{onty7}\frac{1}{\bar{b}}\frac{(1-\bar{b})(r^{1+\alpha_1}-1)}{(r-1)(1-(1-\bar{b})r^{1+\alpha_1})}\frac{r(K-\lambda)}{1-\lambda},
\end{eqnarray}
for all $a\in C$, $1\ls r\ls \min\{\lambda^{-1},(1-\bar{b})^{-\frac{1}{1+\alpha_1}}\}$.
\end{pros}
\begin{dwd}
By (\ref{onty2}) we have that 
$$
\bar{H}_1(r,x,0)\ls H_1(r,x,0)+\frac{(1-\bar{b})H_1(r)G(r,x,0)}{1-(1-\bar{b})G(r)}.
$$
Together with $H_1(r,x,0)=H_1(r,x)$ and $G(r,x,0)=G(r,x)\ls V(x)$ for $x\not\in C$ it follows that
$$
\bar{H}_1(r,x,0)\ls H_1(r,x)+\frac{(1-\bar{b})H_1(r)V(x)}{1-(1-\bar{b})G(r)}.
$$
Consequently by Proposition \ref{pro3}
$$
\bar{H}_1(r,x,0)\ls \frac{r\lambda(V(x)-1)}{1-\lambda}+\frac{(1-\bar{b})(G(r)-1)rV(x)}{(r-1)(1-(1-\bar{b})G(r))}.
$$
Using $G(r)\ls r^{1+\bar{\alpha}}$ we deduce that 
\be\label{elf1}
\bar{H}_1(r,x,0)\ls 
\frac{r\lambda(V(x)-1)}{1-\lambda}+
\frac{(1-\bar{b})(r^{1+\alpha_1}-1)rV(x)}{(r-1)(1-(1-\bar{b})r^{1+\alpha_1})},
\ee
for all $x\not\in C$. On the other hand if
$x\in C$, then (\ref{onty2}) implies that
$$
\bar{H}_1(r,x,0)\ls H_1(r)(1+\frac{(1-\bar{b})G(r)}{1-(1-\bar{b})G(r)})=
\frac{H_1(r,x,0)}{1-(1-\bar{b})G(r)},\;\; \mbox{for}\;\;x\in C,
$$
Hence again by $G(r)\ls r^{1+\alpha_1}$
\be\label{elf2}
(1-\bar{b})\bar{H}_1(r,x,0)\ls \frac{(1-\bar{b})(r^{1+\alpha_1}-1)r}{(r-1)(1-(1-\bar{b})r^{1+\alpha_1})}
\mbox{for}\;x\in C.
\ee
From (\ref{elf1}) and (\ref{elf2}) we conclude (\ref{onty5}).
\smallskip

\noindent
Observe that (\ref{onty1}) and (\ref{superlaska}) imply that
\be\label{superlaska2}
\bar{H}_1(r,a,1)=\frac{r(\bar{G}(r,a,1)-1}{r-1}\ls \frac{\sup_{x\in C}r(G(r,x)-1)}{(r-1)(1-(1-\bar{b})G(r))}=\frac{\sup_{x\in C}H_1(r,x)}{1-(1-\bar{b})G(r)}
\ee
and therefore
$$
\bar{H}_1(r,a,1)\ls \frac{\sup_{x\in C}(H_1(r,x))}{1-(1-\bar{b})r^{1+\alpha_1}}\ls\frac{r(K-\lambda)}{(1-\lambda)(1-(1-\bar{b})r^{1+\alpha_1})},
$$
which is (\ref{onty6}). To prove the last assertion we 
use (\ref{onty21}) and (\ref{superlaska2}), which imply that 
\begin{align*}
& \bar{H}_1(r,a,1)-r\bar{H}_1(1,a,1)\ls 
\frac{1}{\bar{b}}\sup_{x\in C}(H_1(r,x)-rH_1(r,x))+\\
&+\frac{(1-\bar{b})H_1(r)\sup_{x\in C}(G(r,x)-1)}{1-(1-\bar{b})G(r)}.
\end{align*}
The above inequality is equivalent to
\begin{align*}
&\bar{H}_1(r,a,1)-r\bar{H}_1(1,a,1)\ls \frac{1}{\bar{b}}\sup_{x\in C}(H_1(r,x)-rH_1(r,x))+\\
&+\frac{(1-\bar{b})(G(r)-1)\sup_{x\in C}H_1(r,x)}{1-(1-\bar{b})G(r)}.
\end{align*}
Clearly $H_1(r,x)=\frac{r(G(r,x)-1)}{r-1}$, thus by Proposition \ref{pro3} we obtain that
\begin{align*}
& \frac{\bar{H}_1(r,a,1)-r\bar{H}_1(1,a,1)}{r-1}\ls \\
&\ls \frac{1}{\bar{b}}\frac{r\lambda(K-1)}{(1-\lambda)^2}
+\frac{(1-\bar{b})(G(r)-1)}{\bar{b}(r-1)(1-(1-\bar{b})G(r))}\frac{r(K-\lambda)}{1-\lambda}.
\end{align*}
Due to $G(r)\ls r^{1+\alpha_1}$ we deduce (\ref{onty7}) and complete the proof of the result.
\end{dwd}
The consequence of Propositions \ref{pro5}, \ref{pro7} and \ref{pro8} 
is the following theorem in the non-atomic case.
\bt\label{thm3}
Suppose $(X_n)_{n\gs 0}$ satisfies conditions \ref{A1}-(\ref{A3} form the introduction.  Then $(X_n)_{n\gs 0}$
is geometrically ergodic - it verifies (\ref{conv1}) and we have the following bounds 
on $\rho_V$, $M_1$:
$$
\rho_V\ls \bar{r}_0^{-1}
$$ 
and  
\begin{align*}
&M_1(r)\ls \frac{2\lambda r }{1-\lambda}+\frac{2(1-\bar{b})(r^{1+\alpha_1}-1)r}{(r-1)(1-(1-\bar{b})r^{1+\alpha_1})}
+\frac{\bar{b}}{1-(1-\bar{b})r^{1+\alpha_1}}
\frac{r\lambda(K-1)}{(1-\lambda)^2}+\\
& +\frac{(r-1)\bar{K}_0(r)+\bar{b}(1-\bar{b})(r^{1+\alpha_1}-1)}{(r-1)(1-(1-\bar{b})r^{1+\alpha_1})^2}\frac{r(K-\lambda)}{1-\lambda},
\end{align*}
where $\bar{K}_0(r)=\bar{K}_0(r,b,\bar{b},\lambda^{-1},K\lambda^{-1})$, $\bar{r}_0=\bar{r}_0(b,\bar{b},\lambda^{-1},K\lambda^{-1})$ 
are given in Theorems \ref{cori6} and \ref{cori7}.
\et
\begin{dwd}
Note that $\bar{b}\pi(C)\bar{H}_1(1,a,1)=1$.
We apply Proposition \ref{pro5} in the way that
we sum (\ref{bam1}) with weight $(1-\bar{b}1_{C}(x))$
and (\ref{bam2}) with weight $\bar{b}1_{C}(x)$. 
Then we use (\ref{onty4}) to bound 
$\bar{b}1_{C}(x)+(1-\bar{b}1_C(x))\bar{G}(r,x,0)$ and 
(\ref{onty5}) to bound $(1-\bar{b}1_C(x))\bar{H}_1(r,x,0)=(1-\bar{b}1_C(x))\frac{r\bar{G}(r,x,0)-1}{r-1}$.
Finally (\ref{onty6}) and (\ref{onty7}) are
estimates for $\bar{H}(r,a,1)$ and $(\bar{H}_1(r,a,1)-r\bar{H}_1(1,a,1))/(r-1)$.
\end{dwd}

\subsection{Case of $W\equiv V$}

The second case we consider is when $W=V$. 
\begin{pros}\label{pro9}
The following inequalities hold
\be\label{onty8}
(1-\bar{b}1_C(x))\bar{H}_V(r,x,0)\ls 
\frac{\lambda r(V(x)-1)}{1-r\lambda}+
(\frac{K-r\lambda}{1-r\lambda}-\bar{b})\frac{rV(x)}{1-(1-\bar{b})r^{1+\alpha_1}},
\ee
for all $x\in \ccS$, $1\ls r\ls \min\{\lambda^{-1},(1-\bar{b})^{-\frac{1}{1+\alpha_1}}\}$,
\be\label{onty9}
\bar{H}_V(r,a,1)\ls
\bar{b}^{-1}\frac{r(K-r\lambda)}{1-r\lambda}+
\bar{b}^{-1}(\frac{K-r\lambda}{1-r\lambda}-\bar{b})\frac{r-1}{1-(1-\bar{b})r^{1+\alpha_1}}\frac{r(K-\lambda)}{1-\lambda},
\ee
for all $a\in C$, $1\ls r\ls \min\{\lambda^{-1},(1-\bar{b})^{-\frac{1}{1+\alpha_1}}\}$, in particular $\bar{H}_V(1,a,1)\ls \bar{b}^{-1}\frac{K-\lambda}{1-\lambda}$,
and
\begin{eqnarray}
&& \frac{\bar{H}_V(r,a,1)-r\bar{H}(1,a,1)}{r-1}\ls 
\bar{b}^{-1}\frac{r(K-1)}{(1-\lambda)(1-r\lambda)}+\nonumber\\
&&\label{onty10}+\bar{b}^{-1}(\frac{K-r\lambda}{1-r\lambda}-\bar{b})\frac{1}{1-(1-\bar{b})r^{1+\alpha_1}}\frac{r(K-\lambda)}{1-\lambda},
\end{eqnarray}
for all $a\in C$, $1\ls r\ls \min\{\lambda^{-1},(1-\bar{b})^{-\frac{1}{1+\alpha_1}}\}$.
\end{pros}
\begin{dwd}
We recall that (\ref{onty2}) implies that
$$
\bar{H}_V(r,x,0)\ls H_V(r,x,0)+\frac{(1-\bar{b})H_V(r)G(r,x,0)}{1-(1-\bar{b})G(r)}.
$$
Therefore since $H_V(r,x,0)=H_V(r,x)$ and $G(r,x,0)=G(r,x)$ for all $x\not\in C$
we can use Propositions \ref{pro1} and  \ref{pro4} to get
$$
\bar{H}_V(r,x,0)\ls \frac{\lambda r(V(x)-1)}{1-r\lambda}+\frac{(1-\bar{b})H_V(r)V(x)}{1-(1-\bar{b})G(r)}
$$
for all $x\not\in C$. Similarly for $x\in C$
$$
(1-\bar{b})\bar{H}_V(r,x,0)\ls (1-\bar{b})H_V(r)(1+\frac{(1-\bar{b})G(r)}{1-(1-\bar{b})G(r)})=\frac{(1-\bar{b})H_V(r)}{1-(1-\bar{b})G(r)}.
$$
Hence using $G(r)\ls r^{1+\alpha_1}$ we obtain that
$$
(1-\bar{b}1_{C}(x))\bar{H}_V(r,x,0)\ls \frac{\lambda r(V(x)-1)}{1-r\lambda}+\frac{(1-\bar{b})H_V(r)V(x)}{1-(1-\bar{b})r^{1+\alpha_1}}.
$$
Therefore it suffices to bound $(1-\bar{b})H_V(r)$, note that
$$
\bar{b}H_V(r,x,1)+(1-\bar{b})H_V(r)\ls \sup_{x\in C}H_V(r,x),\;\;\mbox{for all}\;x\in C.
$$
Clearly $H_V(r,x,1)\gs r$, so by Proposition \ref{pro4} we deduce that
\be\label{iji}
(1-\bar{b})H_V(r)\ls \frac{r(K-r\lambda)}{1-r\lambda}-\bar{b}r,
\ee
which establishes (\ref{onty8}). To show the remaining assertions we use bounds (\ref{onty11}), (\ref{onty21})
and (\ref{superlaska}) obtaining that
$$
\bar{H}_V(r,a,1)\ls \bar{b}^{-1}\sup_{x\in C}H_V(r,x)+\bar{b}^{-1}\frac{(1-\bar{b})H_V(r)\sup_{x\in C}(G(r,x)-1)}{1-(1-\bar{b})G(r)} 
$$
and
\begin{align*}
&\bar{H}_V(r,a,1)-r\bar{H}_V(r,a,1)\ls \bar{b}^{-1}\sup_{x\in C}(H_V(r,x)-rH_V(1,x))+\\
&+\bar{b}^{-1}\frac{(1-\bar{b})H_V(r)\sup_{x\in C}(G(r,x)-1)}{1-(1-\bar{b})G(r)} 
\end{align*}
Recall that $G(r)\ls r^{1+\alpha_1}$ and by Propositions \ref{pro1} and \ref{pro4}
$$
\frac{G(r,x)-1}{r-1}\ls \frac{K-\lambda}{1-\lambda},\;\;
H_V(r,x)\ls \frac{r(K-r\lambda)}{1-r\lambda},
$$
consequently
$$
\bar{H}_V(r,a,1)\ls \bar{b}^{-1}\frac{r(K-r\lambda)}{1-r\lambda}+
\bar{b}^{-1}\frac{(r-1)(1-\bar{b})H_V(r)}{1-(1-\bar{b})r^{1+\alpha_1}}\frac{K-\lambda}{1-\lambda}. 
$$
Together with (\ref{iji}) it completes the proof of (\ref{onty9}). Finally the same argument shows
$$
\frac{\bar{H}_V(r,a,1)-r\bar{H}_V(r,a,1)}{r-1}\ls \bar{b}^{-1}\frac{r(K-1)}{(1-\lambda)(1-r\lambda)}+
\bar{b}^{-1}\frac{(1-\bar{b})H_V(r)}{1-(1-\bar{b})r^{1+\alpha_1}}\frac{K-\lambda}{1-\lambda}.
$$
Again by (\ref{iji}) we obtain (\ref{onty10}), which completes the proof.
\end{dwd}
Propositions \ref{pro5}, \ref{pro7} and \ref{pro9} imply the estimate in the most general case 
when there is no control on $W$ but $V$.
\bt\label{thm4}
Suppose $(X_n)_{n\gs 0}$ satisfies (\ref{A1})-(\ref{A3}).  Then $(X_n)_{n\gs 0}$
is geometrically ergodic - it verifies (\ref{conv1}) and we have the following bounds 
on $\rho_V$, $M_1$:
$$
\rho_V\ls \bar{r}_0^{-1}
$$ 
and 
\begin{align*}
& M_V(r)\ls
\frac{\lambda r}{1-r\lambda}+(\frac{K-r\lambda}{1-r\lambda}-\bar{b})\frac{r}{1-(1-\bar{b})r^{1+\alpha_1}}+\\
&+\frac{K-\lambda}{1-\lambda}(\frac{r\lambda}{1-\lambda}+\frac{(1-\bar{b})(r^{1+\alpha_1}-1)r}{(r-1)(1-(1-\bar{b})r^{1+\alpha_1})})+\\
&+\frac{\bar{b}}{1-(1-\bar{b})r^{1+\alpha_1}}(\frac{r(K-1)}{(1-\lambda)(1-r\lambda)}+\\
&+(\frac{K-r\lambda}{1-r\lambda}-\bar{b})\frac{1}{1-(1-\bar{b})r^{1+\alpha_1}}\frac{r(K-\lambda)}{1-\lambda})+\\
&+\frac{\bar{K}_0(r)}{1-(1-\bar{b})r^{1+\alpha_1}}
(\frac{r(K-r\lambda)}{1-r\lambda}+(\frac{K-r\lambda}{1-r\lambda}-\bar{b})\frac{r-1}{1-(1-\bar{b})r^{1+\alpha_1}}\frac{r(K-\lambda)}{1-\lambda}),
\end{align*}
where $\bar{K}_0(r)=\bar{K}_0(r,b,\bar{b},\lambda^{-1},K\lambda^{-1})$, $\bar{r}_0=\bar{r}_0(b,\bar{b},\lambda^{-1},K\lambda^{-1})$ 
are given in Theorems \ref{cori6} and \ref{cori7}.
\et 
\begin{dwd}
Observe that $\pi(C)\ls 1$. As in the proof of Theorem \ref{thm3} we 
we use Proposition \ref{pro5} summing (\ref{bam1}) with weight $(1-\bar{b}1_{C}(x))$
and (\ref{bam2}) with weight $\bar{b}1_{C}(x)$. 
Again we use (\ref{onty4}) to bound 
$\bar{b}1_{C}(x)+(1-\bar{b}1_C(x))\bar{G}(r,x,0)$, then
(\ref{onty8}), (\ref{onty9}), (\ref{onty10}) to bound respectively $(1-\bar{b}1_C(x))\bar{H}_V(r,x,0)$,
$\bar{H}_V(r,a,1)$ and $(\bar{H}_V(r,a,1)-r\bar{H}_V(1,a,1))/(r-1)$. We use also the bound
$\bar{H}_V(1,a,1)\ls \bar{b}^{-1}\frac{K-\lambda}{1-\lambda}$ and 
(\ref{onty5}) to bound $(1-\bar{b}1_C(x))\bar{H}_1(r,x,0)=(1-\bar{b}1_C(x))\frac{r\bar{G}(r,x,0)-1}{r-1}$.
\end{dwd}

\section*{Appendix B}\label{appendB}

\renewcommand{\thethm}{B.\arabic{thm}}
\setcounter{thm}{0}
\renewcommand{\thesubsection}{B.\arabic{subsection}}
\setcounter{subsection}{0}
\renewcommand{\theequation}{B.\arabic{equation}}
\setcounter{equation}{0}

We compare our result with what was shown in \cite{Bax} as a numerical
test for the presented approach.
\subsection{\textit{The reflecting Random Walk}}
We consider the Bernoulli random walk on $\Z_{+}$ with transition probabilities $P(i,i-1)=p>1/2$, $P(i,i+1)=q=1-p$ for
$i\gs 1$ and boundary conditions $P(0,0)=p$, $P(0,1)=q$. We set $C=\{0\}$ and $V(i)=(p/q)^{i/2}$, and compute
$\lambda=2\sqrt{pq}$, $K=p+\sqrt{pq}$, $b=p$ and $u_{\infty}=\pi(C)=1-q/p$. The optimal radius of convergence
for the reflecting random walk is $\lambda$. 
\smallskip

\noindent
Consider two cases:
\begin{enumerate}
\item In the first one we consider
$p=2/3$, so $b=2/3$, $\lambda=2\sqrt{2}/3$, $K=(2+\sqrt{2})/3$, $u_{\infty}=1/2$. 
\item In the second one we set
$p=0.9$, and hence $\lambda=0.6$, $K=1.2$, $b=0.9$, $u_{\infty}=8/9$.
\end{enumerate}
We compare our result in Table 1 below, where $\rho$, $\rho_C$ denotes
estimates on the radius of convergence
in the case where respectively $u_{\infty}$ is known or not. We use \textit{Optimal} for the true value of the spectral radius and \textit{Bednorz}, \textit{Baxendale},
\textit{Meyn-Tweedie1} and \textit{Meyn-Tweedie2} respectively for our Corollaries \ref{cori1} and \ref{cori2}, Baxendale's Theorem 3.2 \cite{Bax},
Meyn-Tweedie's result given in  \cite{MT2} and 
its improved version (see Section 8 in \cite{Bax} for details). 
\begin{align*}
& \mbox{Table 1}  \\
&\begin{array}{llllll}
p=2/3 & \rho & \rho_C \\  
Optimal & 0,9428 & 0,9428 \\
Bednorz & 0,9737 & 0,9737\\
Baxendale & 0,9994 & ? \\
Meyn-Tweedie1 & 0,9999 & 0,9988 \\
Meyn-Tweedie2 & 0,9991 & 0,9927
\end{array},\\
&\begin{array}{llllll}
p=0,9 & \rho & \rho_C \\  
Optimal &  0,6 & 0,6 \\
Bednorz & 0,6 & 0,6 \\
Baxendale & 0,9060 & ?\\
Meyn-Tweedie1 & 0,9967 & 0,9888\\
Meyn-Tweedie2 & 0,9470 & 0,9467
\end{array} 
\end{align*}

\subsection{Metropolis Hastings Algorithm for the Normal Distribution}

In this example we consider the convergence of a Metropolis-Hastings algorithm in the case 
when we want to simulate $\pi=\ccN(0,1)$ with candidate transition probability $q(x,\cdot)=\ccN(x,1)$.
The example was studied in \cite{MT2} and also in \cite{Rob1} and \cite{Rob2}.
By the algorithm definition $\P(x,\cdot)$ is distributed with a density
$$
p(x,y)=\begin{array}{ll}
\frac{1}{\sqrt{2\pi}}\exp(-\frac{(y-x)^2}{2}), &  \mbox{if}\;\;|x|\gs |y|\\
\frac{1}{\sqrt{2\pi}}\exp(-\frac{(y-x)^2+y^2-x^2}{2}), & \mbox{if}\;\; |x|\ls |y|.
\end{array}
$$
The natural setting of the problem is 
to consider Lyapunov functions of the type $V(x)=e^{s|x|}$ and $C=[-d,d]$.
Consequently (see \cite{Bax} for details)
$$
\lambda=\frac{PV(d)}{V(d)},\;\; K=PV(d)=e^{sd}\lambda.
$$
The computed value for $\rho$ depends on $d$ and $s$, and hence to 
we need to find the optimal ones. Moreover to compare our result with the previous
contributions to the problem, let $\nu$ be given by
$$
\nu(dx)=c\exp(-x^2)1_{C}(x)dx,
$$
for a suitable normalizing constant $c$. In this case, $\nu(C)=1$ and we have 
$$
b=\bar{b}=\sqrt{2}\exp(-d^2)[\Phi(\sqrt{2}d)-1/2].
$$ 
In this case we work with  the additional complication of the splitting construction. The results are compared in Table 2, where again \textit{Bednorz1} and
\textit{Bednorz2} denote our Theorems \ref{cori7} and \ref{cori6}
(depending whether or not we use the additional information on $\pi(C)$), 
\textit{Baxendale} denotes what can be obtained by Baxendale's
Theorem 3.2,\textit{Coupling} denotes the estimate obtained by the coupling approach (see Section 7 in \cite{Bax} and \cite{Ros2}) and \textit{Meyn-Tweedie} the result obtained in the original paper \cite{MT2}. Note that we compare methods that no additional assumptions
on the transition probabilities are made. 
\begin{align*}
& \mbox{Table 2}  \\
&\begin{array}{llll}
 & d & s & (1-\rho) \\  
Bednorz1 & 0.96 & 0.065 &     0.00000529\\
Bednorz2 & 0.92 & 0.169 &  0.00005496 \\
Baxendale & 1 & 0.13 & 0.00000063\\
Coupling & 1.8 & 1.1 & 0.00068\\
Meyn-Tweedie & 1.4 & 0.00004 & 0.000000016
\end{array}.
\end{align*}
Another possible choice of $\nu$ is
$$
\bar{b}\nu(dx)=\left\{
\begin{array}{ll} \frac{1}{\sqrt{2\pi}}\exp(-\frac{(|x|+d)^2}{2})dx, & \mbox{if}\;\; |x|\ls d \\
\frac{1}{\sqrt{2\pi}}\exp(-d|x|-|x|^2)dx, &
\mbox{if}\;\; |x|\gs d \end{array}\right.
$$
In this case
$$
b=2(\Phi(2d)-\Phi(d))\;\;\mbox{and}\;\;
\bar{b}=b+\sqrt{2}\exp(d^2/4)(1-\Phi(3d/\sqrt{2})).
$$
Using the same notation as in Table 2 we compare the results below.
\begin{align*}
& \mbox{Table 3}  \\
&\begin{array}{llll}
 & d & s & (1-\rho) \\  
Bednorz1 & 1.03 & 0,0733 & 0,00001061 \\
Bednorz2 & 0.97 & 0,1740 & 0,00013637 \\
Baxendale & 1 & 0,16 & 0,0000017\\
Coupling & 1,9 & 1,1 & 0,00187\\
\end{array}.
\end{align*}
Observe that our method is bit worse than coupling yet
it is relatively simple (does not require further examination of the Lyapunov function $V$.)

\subsection{Contracting Normals}

Here we consider the family of Markov chains with transition probability $\P(x,\cdot)=\ccN(\theta x,1-\theta^2)$ for some parameter $\theta\in (-1,1)$.
The family of examples occurs in \cite{Ros} as one 
component of a two component Gibbs sampler.
The example was discussed in \cite{Bax}, \cite{Rob1}
and \cite{Rob2}.
Here we take
$V(x)=1+x^2$ and $C=[-c,c]$. Then (\ref{A2}) is satisfied
with
$$
\lambda=\theta^2+2\frac{1-\theta^2}{1+c^2},\;\;
K=2+\theta^2(c^2-1).
$$
The we choose $\nu$ concentrated on $C$ so that
$$
\bar{b}\nu(dy)=\min_{x\in C}\frac{1}{\sqrt{2\pi(1-\theta^2)}}\exp(-\frac{(\theta x-y)^2}{2(1-\theta^2)})dy 
$$
for $y\in C$. Integrating with respect to $y$ gives
$$
\bar{b}=2(\Phi(\frac{(1+|\theta|)c}{\sqrt{1-\theta^2}})-
\Phi(\frac{|\theta|c}{\sqrt{1-\theta^2}})).
$$
We compare our answer \textit{Bednorz1}, \textit{Bednorz2} (Theorems \ref{cori7}, \ref{cori6} resp.) with the coupling method \textit{Coupling} and \textit{Baxendale2} an approach based on Kendal type result (Theorem 3.3 in \cite{Bax}) that requires invertibility of the transition function.
\begin{align*}
& \mbox{Table 4}  \\
&\begin{array}{llll}
 & \theta & c & 1-\rho \\ 
Bednorz1 & 0,5 & 1,5 &  0,000872023152 \\
Bednorz1 & 0,75 & 1,2 & 0,000000964524 \\
Bednorz1 & 0,9 & 1,1 &  0,000000000004 
\end{array},\\
& \begin{array}{llll}
 & \theta & c & 1-\rho \\  
Bednorz2 & 0,5 & 1,5 &  0,002754672439 \\
Bednorz2 & 0,75 & 1,2 & 0,000017954821 \\
Bednorz2 & 0,9 & 1,1 &  0,000000000881
\end{array}, \\
& \begin{array}{llll}
 & \theta & c & 1-\rho \\  
Baxendale2 & 0,5 & 1,5 & 0,050 \\
Baxendale2 & 0,75 & 1,2 & 0,0042 \\
Baxendale2 & 0,9 & 1,1 & 0,00002\\
\end{array},\\
& \begin{array}{llll}
 & \theta & c & 1-\rho \\  
Coupling & 0,5 & 2,1 & 0,054 \\
Coupling & 0,75 & 1,7 & 0,0027 \\
Coupling & 0,9 & 1,5 & 0,00002\\
\end{array}
\end{align*}

\subsection{Reflecting random walk, continued.}

Here we slightly redefine our first example. Let
$\P(0,\{0\})=1$ and $\P(0,\{1\})=1-\va$ for some $\va>0$. We concentrate on the difficult case, 
when  $\va<p$, that was studied in \cite{Rob2} and \cite{For}. Note that 
when $\va\gs p$ then the chain is stochastically monotone and then the result of Tweedie \cite{LT} apply. Let $V(i)=(p/q)^{i/2}$, $C=\{0\}$ as earlier.
Then $\lambda=2\sqrt{pq}$, $K=\va+(1-\va)\sqrt{p/q}$ and
$b=\va$. In this example we can calculate the formula on $b(z)$ which is
\be\label{miz}
b(z)=G(z,0)=\va z+(1-\va)zG(z,1)=\va z+\frac{1-\va}{2q}(1-(1-4pqz^2)^{1/2}),
\ee
for $|z|<1/\sqrt{4pq}$, where the formula for $G(z,1)$ is from \cite{Fel}.
Consequently
$$
\pi(\{0\})^{-1}=b'(1)=\va+\frac{2p(1-\va)}{p-q}.
$$
On the other hand (\ref{miz}) leads to the optimal bound on the radius on convergence
which is
$$
\rho=\left\{\begin{array}{ll}
\frac{pq+(p-\va)^2}{p-\va}, & \mbox{if}\;\;\va<\frac{p-q}{1+\sqrt{q/p}},\\
2\sqrt{pq}, &\mbox{otherwise}.\end{array} \right.
$$
We compare \textit{Bednorz1},\textit{Bednorz} - our Corollaries \ref{cori1},\ref{cori2}
with results \textit{Fort}
and \textit{Baxendale}
that denotes respectively 
the result of Fort \cite{For} and Baxendale's Theorem 1.2
\cite{Bax}. Note that both two methods use further properties of transition probability in this particular example.
\begin{align*}
& \mbox{Table 5} \\
& \begin{array}{lllll}
p=0,6 & \va=0,05 & \va=0,25 & \va=0,5  \\  
Optimal & 0,9864 & 0,9798 & 0,9798  \\
Bednorz1 & 0,99993 & 0,9994 & 0,99783\\
Bednorz2 & 0,99993 & 0,9994 & 0,9977\\
Fort & 0,9997 & 0,9995 & 0,9994\\
Bax & 0,9909 & 0,9798 & 0,9798 
\end{array},\\
&\begin{array}{lllll}
p=0,7 & \va=0,05 & \va=0,25 & \va=0,5 \\  
Optimal & 0,9165 & 0,9165 & 0,9165 \\
Bednorz1 & 0,9992 & 0,9940 & 0.9783\\
Bednorz2 & 0,9991 & 0,9935 & 0.9779 \\
Fort & 0,9964 & 0,9830 & 0,9757 \\
Bax & 0,9731 & 0,9165 & 0,9165 
\end{array},\\ 
& \begin{array}{lllll}
p=0,8 & \va=0,05 & \va=0,25 & \va=0,5 \\  
Optimal & 0,9633 & 0,8409 & 0,8000 \\
Bednorz1 & 0,9970 & 0,9780 & 0,9266\\
Bednorz2 & 0,9964 & 0,9751 & 0,9253 \\
Fort & 0,9793 & 0,9333 & 0,9333 \\
Bax & 0,9759 & 0,8796 & 0,8000  
\end{array},\\
& \begin{array}{lllll}
p=0,9 & \va=0,05 & \va=0,25 & \va=0,5 \\  
Optimal & 0,9559 & 0,7885 & 0,6250 \\
Bednorz1 & 0,9927 & 0,9489 & 0,8408\\
Bednorz2 & 0,9899 & 0,9358 & 0,8280\\
Fort & 0,9696 & 0,8539 & 0,7500 \\
Bax & 0,9687 & 0,8470 & 0,6817
\end{array},   \\
& \begin{array}{lllll}
p=0,95 & \va=0,5 & \va=0,25 & \va=0,5 \\  
Optimal & 0,9528 & 0,7679 & 0,5556 \\
Bednorz1 & 0,9888 & 0,9249 & 0,7827 \\
Bednorz2 & 0,9841 & 0,9024 &  0,7537 \\
Fort & 0,9564 & 0,7853 & 0,5814 \\
Bax & 0,9645 & 0,7853 & 0,5814  
\end{array}.  
\end{align*}

\end{document}